\documentclass{amsart}

\usepackage{amsmath,amsthm}
\usepackage{latexsym}
\usepackage{amssymb}
\newcommand{\less}{\lesssim}

\newcommand{\beal}{\begin{align}}
\newcommand{\enal}{\end{align}}
\newcommand{\bealn}{\begin{align*}}
\newcommand{\enaln}{\end{align*}}
\newcommand{\bear}{\begin{eqnarray}}
\newcommand{\eear}{\end{eqnarray}}
\newcommand{\beeq}{\begin{equation}}
\newcommand{\eneq}{\end{equation}}

\newcommand{\eps}{{\varepsilon}}
\newcommand{\R}{{\mathbb R}}

\newcommand{\Z}{{\mathbb Z}}
\newcommand{\Nat}{{\mathbb N}}
\newcommand{\Compl}{{\mathbb C}}
\newcommand{\calM}{{\mathcal M}}
\newcommand{\la}{\langle}
\newcommand{\ra}{\rangle}

\newcommand{\etao}{\xi_d}
\newcommand{\go}{\phi_d}
\def\pr{\partial}

\def\nn{\nonumber}
\def\calD{\mathcal D}
\def\calH{\mathcal H}
\def\f{\frac}
\def\lam{\lambda}

\def\bvec{\left[ \begin{array}{c}}
\def\endvec{\end{array}\right]}
\def\spec{{\rm spec}}

\def\calF{{\mathcal F}}

\def\Im{\;{\rm Im}\,}
\def\Re{\;{\rm Re}\,}

\def\cQ{{\mathcal Q}}
\def\calE{{\mathcal E}}
\def\calC{{\mathcal C}}

\newcommand{\IS}{\mbox{\,{\rm IS}}}

\def\einbet{\hookrightarrow}
\def\pr{\partial}

\def\nn{\nonumber}
\def\bm{\left[ \begin{array}{cc}}
\def\endm{\end{array}\right]}

\def\f{\frac}
\def\drb{\beta}
\def\drnu{\beta_0}
\renewcommand{\ln}{\log}
\def\cL{{\mathcal L}}
\def\tileps{\tilde\eps}
\def\dom{{\rm Dom}}
\def\cK{{\mathcal K}}
\def\calK{\cK}

\def\Bla{\Big\la}
\def\Bra{\Big\ra}
\def\tilphi{\phi}
\def\calS{{\mathcal S}}
\def\dist{{\rm dist}}

\newtheorem{theorem}{Theorem}
\newtheorem{lemma}[theorem]{Lemma}
\newtheorem{defi}[theorem]{Definition}
\newtheorem{cor}[theorem]{Corollary}
\newtheorem{prop}[theorem]{Proposition}
\newtheorem{proposition}[theorem]{Proposition}

\theoremstyle{remark}
\newtheorem{remark}[theorem]{Remark}

\numberwithin{equation}{section}
\numberwithin{theorem}{section}

\begin{document}

\title[Slow blow-up solutions]{Slow blow-up solutions for the $H^1(\R^3)$ critical
focusing semi-linear wave equation in $\R^3$}

\author{J.\ Krieger}
\address{Harvard University, Dept. of Mathematics, Science Center, 1 Oxford Street, Cambridge, MA 02138, U.S.A.}
\email{jkrieger@math.harvard.edu}

\author{W.\ Schlag}
\address{Department of Mathematics, The University of Chicago, 5734 South University Avenue, Chicago, IL 60637, U.S.A.}
  \email{schlag@math.uchicago.edu}

\author{D.\ Tataru}
\address{Department of Mathematics, The University of California at Berkeley, Evans Hall, Berkeley, CA 94720, U.S.A.}
\email{tataru@math.berkeley.edu}

\thanks{The authors were partially supported by the National Science Foundation, J.\ K.\ by DMS-0401177, W.\ S.\ by
DMS-0617854, D.~T.\ by DMS-0354539, and DMS-0301122.}

\maketitle

\begin{abstract} Given $\nu>\frac12$ and $\delta>0$ arbitrary,
we prove the existence of energy solutions of
\begin{equation}\label{eq:wave_abstr}
\partial_{tt} u - \Delta u - u^5 =0
\end{equation}
in $\R^{3+1}$ that blow up exactly at $r=t=0$ as $t\to0-$. These
solutions are radial and of the form $u = \lambda(t)^{\frac12}
W(\lambda(t)r) + \eta(r,t)$ inside the cone $r\le t$, where
$\lambda(t)=t^{-1-\nu}$, $W(r)=(1+r^2/3)^{-\frac12}$ is the
stationary solution of~\eqref{eq:wave_abstr}, and $\eta$ is a
radiation term with
\[
\int_{[r\le t]} \big(|\nabla \eta(x,t)|^2 +
|\eta_t(x,t)|^2+|\eta(x,t)|^6\big)\, dx \to 0, \qquad t\to0
\]
Outside of the light-cone there is the energy bound
\[
\int_{[r>t]} \big( |\nabla u(x,t)|^2+|u_t(x,t)|^2+|u(x,t)|^6\big)\,
dx <\delta
\]
for all small $t>0$. The regularity of $u$ increases with $\nu$. As
in our accompanying paper on wave-maps~\cite{KST}, the argument is
based on a renormalization method for the `soliton profile' $W(r)$.
\end{abstract}

\section{Introduction}

Since the seminal paper of J\"orgens \cite{Jor} much work has been
devoted to the study of well-posedness of the nonlinear wave
equation
\[
\partial_{tt} u -
\Delta u + f(u)=0
\]
in $\R^{3+1}_{x,t}$ and suitable nonlinearities $f(u)$. J\"orgens
showed that for $H^1(\R^3)$ subcritical defocusing nonlinearities
$f(u)=|u|^{p-1} u$ with $p<5$ smooth data lead to smooth solutions
for all times. The critical defocusing case $p=5$ was resolved by
Struwe~\cite{Struwe} for radial data and Grillakis~\cite{Gril} for
general data. These authors proved global well-posedness and
scattering results for energy solutions, see
Shatah--Struwe~\cite{SS} and Sogge~\cite{Sogge}.  No corresponding
results are known for the supercritical case $p>5$.

In this paper we address the solvability of the nonlinear wave
equation in $\R^{3+1}$ with a focusing nonlinearity
$f(u)=-|u|^{p-1}u$. In this case blow-up may occur. Indeed, it was
shown by Levine~\cite{Lev} via a convexity argument that data in
$(\dot{H^1}\cap L^2)\times L^2$ with negative energy lead to
finite-time blow-up, see also Strauss~\cite{Strauss}. Local
well-posedness in the optimal regularity class was considered by
several authors, see Sogge~\cite{Sogge} for a detailed exposition of
this work. Most relevant for us is the case $p=5$ where the equation
is locally well-posed in the energy space $\dot{H}^1\times
L^2(\R^3)$. Moreover, if the solution cannot be continued beyond
some finite time~$T_*$ as an energy solution, then necessarily the
Strichartz norm $\|u\|_{L^8([0,T_*)\times \R^3)} =\infty$ (with
similar results in all dimensions).

The question of the blow-up rate was addressed by Merle--Zaag in the
conformal range $p\le 3$, see~\cite{MZ1}--\cite{MZ3} (their results
extend to all dimensions). They showed that if solutions to the
Cauchy problem, with $1<p\le 3$,
\[
\partial_{tt} u -
\Delta u - |u|^{p-1}u =0,\quad (u_0,u_1)\in H^1_{\rm loc}\times
L^2_{\rm loc}
\]
blow up in finite time $T_*$, then the following holds: for any
$a\in\R^3$ the self-similar change of variables
\[u(x,t)=(T_*-t)^{-\frac{2}{p-1}} w_a(y,s),\quad y=\frac{x-a}{T_*-t},
\; s=-\log(T_*-t)
\]
leads to functions $w_a$ satisfying
\[
\sup_{s\ge -\log T_*+1,\;a\in\R^3} \|w_a(s)\|_{H^1(B)} +
\|\partial_s w_a(s)\|_{L^2(B)} \le K
\]
where $B$ is the unit ball and a constant $K$ that only depends on
$p,T_*$ and the norm of the initial data in $H^1_{\rm loc}\times
L^2_{\rm loc}$.

 For the energy critical case $p=5$, i.e.,
\begin{equation}\label{TheEquation} \partial_{tt} u - \Delta u -
u^5=0 \end{equation}
 there has been some recent activity, see
\cite{KS}, \cite{KM}, \cite{KarStr},  which we now describe in more
detail. The Talenti--Aubin solutions
\[
W(r)=(1+r^2/3)^{-\frac12}
\]
are extremizers of the Sobolev imbedding $\dot{H}^1(\R^3)\einbet
L^6(\R^3)$ and satisfy the Euler-Lagrange equation $-\Delta W -
W^5=0$. In~\cite{KS} the first two authors showed that there exists a
small co-dimension one manifold $\calM$ around $W$ in a suitable
topology so that data on this manifold exhibit global existence and an
asymptotic behavior of bulk term plus radiation. The radiation term is
also shown to scatter like a free energy wave. It is conjectured,
see~\cite{Biz}, that this manifold has the property that it separates
a region of scattering from one of blow-up. As a first result in this
direction, Karageorgis--Strauss~\cite{KarStr} showed that above the
tangent space of $\calM$ at $W$ finite time blow-up occurs, albeit for
the equation
\[
\partial_{tt} u - \Delta u -
|u|^5=0
\]
Note that the result of \cite{KS} equally well applies to the
nonlinearity $|u|^5$ (in fact, the solutions constructed
in~\cite{KS} are positive so that there is no distinction between
$u^5$ and $|u|^5$ from the point of view of that paper).

 Kenig--Merle~\cite{KM} studied the behavior of
solutions with data $(u_0,u_1)\in \dot{H}^1\times L^2(\R^3)$ of
energy $\calE(u_0,u_1)< \calE(W,0)$ where the conserved energy is
\[
\calE(u,u_t) = \int_{\R^3} \Big[\frac12 ( u_t^2 + |\nabla u|^2) -
\frac{|u|^6}{6}\Big] \, dx
\]
They found that in this regime there is a dichotomy between blow-up
and global existence/scattering depending on whether $\|\nabla
u_0\|_2
> \|\nabla W\|_2$ or $\|\nabla u_0\|_2 < \|\nabla W\|_2$.

Note that
\[
W(x,\lambda):= \lambda^{\frac12} W(\lambda x)
\]
is a stationary solution of \eqref{TheEquation} for all $\lambda>0$.
Moreover, the energy is constant in $\lambda$ (reflecting the energy
criticality of the equation).
 Linearizing the wave equation around $W$
leads to the linearized operator
\[ H=-\Delta - 5W^4\]
The wave evolution of $H$ has two types of instabilities: an
 exponential instability arising from the negative spectrum of $H$
(which has a unique negative eigenvalue) as well as  a "bound state"
at zero energy: $H(\partial_\lambda W|_{\lambda=1})=0$ where
$\partial_\lambda W$ decays like $r^{-1}$ and thus does not belong
to $L^2(\R^3)$ --- this is what one refers to as a zero energy
resonance. In this paper we construct blow-up solutions by
`projecting out' the exponentially growing mode of the linearized
equation.

More precisely,  we seek radial, real-valued,  blow-up solutions
\[u(x,t)=\lambda(t)^{\frac12}W(\lambda(t)x)+\eta(x,t)\]
of \eqref{TheEquation}
 where
$\lambda(t)\to\infty$ as $t\to 0$ and with the local energy inside
the light-cone $|x|\le t$ of $\eta(x,t)$ going to zero as $t\to0$.
The {\em local}  energy relative to the origin is defined as
\[
\calE_{\rm loc}(\eta) = \int_{[|x|<t]} ( \eta_t^2 + |\nabla
\eta|^2+|\eta|^6)\, dx
\]
The following theorem is the main result of this paper. The blow-up
occurs at time $t=0$ when solving backwards in time.

\begin{theorem}\label{Main}
Let $\nu>\frac12$ and $\delta>0$. Then there exists an energy
solution $u$ of \eqref{TheEquation} which blows up precisely at
$r=t=0$  and which has the following property: in the cone $|x|=r\le
t $ and for small times $t$ the solution has the form, with
$\lambda(t)=t^{-1-\nu}$,
\[
u(x,t) = \lambda^\frac12(t) W(\lambda(t) r) + \eta(x,t)
\]
where $\calE_{\rm loc}(\eta(\cdot,t))\to0$ as $t\to0$ and outside
the cone $u(x,t)$ satisfies
\[
\int_{[|x|\ge t]} \big[|\nabla u(x,t)|^2+ |u_t(x,t)|^2 +
|u(x,t)|^6\big]\, dx <\delta
\]
for all sufficiently small $t>0$. In particular, the energy of these
blow-up solutions can be chosen arbitrarily close to~$\calE(W,0)$,
i.e., the energy of the stationary solution.
\end{theorem}

The restriction $\nu > 1/2$
  arises only due to technical reasons, and we hope to eliminate it in
  subsequent work. If $\nu>1$, then the solutions from Theorem~\ref{Main}
  belong to $L^\infty(\R^3)$ for all $t>0$ and blow up
at the rate
\[ \|u(\cdot,t)\|_\infty \asymp t^{-(1+\nu)/2}\] as $t\to 0$.
The proof is based on a renormalization procedure analogous to the
one that the authors used for the construction of blow-up solutions
for wave maps in~\cite{KST}. For our purposes this refers to the
fact that we do not simply perturb around $\lambda^\frac12(t)
W(\lambda(t) r)$ to obtain the linearized equation for $\eta$, but
rather first modify the blow-up profile and then perturb around this
"renormalized" profile. More precisely,  fix a large integer $N$.
Then there exists a function $u^e$ satisfying
  \begin{equation}\label{eq:ue}
  u^e\in C^{\frac{\nu+1}{2}-}(\{t_0>t>0, \;|x|\le t\}),\qquad \calE_{\rm
    loc}(u^e)(t)\lesssim (t\lambda(t))^{-1} \text{\  \ as\  \ }t\to0
  \end{equation}
  so that the radiation term $\eta$ above has the form
  \[
  \eta(x,t) = u^e(r,t) + \eps(r,t), \qquad 0\le r\le t
  \]
  where $\eps$ decays at $t=0$. In fact, $\eta$ can be
  extended globally with the property that
\[
\eps \in t^N H^{\frac{\nu+2}{2}-}(\R^3), \qquad  \eps_t \in t^{N-1}
H^{\frac{\nu}{2}-}(\R^3), \qquad \calE_{\rm
    loc}(\eps)(t)\lesssim t^N \text{\  \ as\  \ }t\to0
\]
with spatial norms that are uniformly controlled as $t\to0$.

As this paper was written concurrently with our wave-map
paper~\cite{KST} it is only natural that there would be some
similarities between this paper and~\cite{KST}. In fact, a secondary
goal here is to show that the method used in both papers is flexible
and applies to quite distinct scenarios. The main differences
between this paper and~\cite{KST} are as follows:
\begin{itemize}
\item The blow-up profile is not constant in $L^\infty$ but rather
grows at rate $\lambda^{\frac12}$. The renormalization procedure
thus needs to be adapted to this case.
\item In contrast to~\cite{KST},
the linearized operator exhibits negative spectrum. This produces
exponential instability of the linearized wave flow.
\item The linearized operator no longer exhibits a strongly singular potential
in the sense of~\cite{GZ}. Thus, a (Dirichlet) boundary condition is
needed at $R=0$.
\end{itemize}
We feel that the most important difference listed here is the
exponential instability. In fact, as in the asymptotic stability
paper~\cite{KS}, we need to `project out' this exponential growth.
Our blow-up rates are therefore expected to be non-generic.

\section{The renormalization step}
\label{sec:mod}

In this section we show how to construct an arbitrarily good
approximate radial solution to the wave equation~\eqref{TheEquation}
as a perturbation of a time-dependent ground state profile
\[
u_0 = \lambda^{\frac12}W(R),\quad W(R) = (1+R^2/3)^{-\frac12},\quad
R =r \lambda(t)
\]
with the polynomial timescale
\[
\lambda(t) = t^{-1-\nu}, \qquad \nu>0.
\]

\begin{theorem}
  \label{thm:sec2} Let $k \in \Nat$. There exists an approximate
  solution $u_{2k-1}$ for \eqref{TheEquation} of the form
  \[
  u_{2k-1}(r,t) = \lambda^{\frac12}(t)\Big[W(R) + \frac{c}{(t\lambda)^2} R^2(1+ R^2)^{-1/2} +
  O\left( \frac{R^2 (1+R^2)^{-\frac32}}{(t\lambda)^2}\right)\Big]
  \]
  so that the corresponding error has size
  \[
  e_{2k-1} = O\left(\frac{\lambda^{\frac12}R}{t^2 (t \lambda)^{2k}}\right)
  \]
 Here the $O(\cdot)$ terms are uniform in $0\le r\le t$ and
  $0<t<t_0$ where $t_0$ is a fixed small constant.
\end{theorem}

\begin{remark}
The $u^e$ in \eqref{eq:ue} is
\[
u^e(r,t)= \lambda^{\frac12}(t)\Big[\frac{c}{(t\lambda)^2} R^2(1+
R^2)^{-1/2} +
  O\left( \frac{R^2 (1+R^2)^{-\frac32}}{(t\lambda)^2}\right)\Big]
\]
The analysis below shows that it has the stated regularity up to the
light-cone. Moreover, one checks that
\[
\calE_{\rm loc}\Big(\frac{\lambda^{\frac12}}{(t\lambda)^2} R
\Big)\less (t\lambda)^{-1}
\]
which is the claimed decay rate for the local kinetic energy of
$u^e$. The local {\em potential} energy of $u^e$ decays  like
$(t\lambda)^{-3}$.
\end{remark}

\begin{proof}
  We iteratively construct a sequence $u_k$ of better approximate
  solutions by adding corrections $v_k$,
  \[
  u_k = v_{k} + u_{k-1}
  \]
  The error at step $k$ is
  \[
  e_k = (-\partial_t^2 + \partial_r^2 +\frac{2}r \partial_r) u_k
  + u_k^5
  \]
  If $u$ were an exact solution, then the difference
  \[
  \eps = u-u_{k-1}
  \]
  would solve the equation
  \begin{equation}\label{eq:gleich}
  (-\partial_t^2 + \partial_r^2 +\frac{2}r \partial_r) \eps + 5u_{k-1}^4 \eps +
  10 u_{k-1}^3 \eps^2 + 10 u_{k-1}^2 \eps^3 + 5u_{k-1} \eps^4 + \eps^5 +
  e_{k-1} =0
  \end{equation}
  In a first approximation we linearize this equation around $\eps
  =0$ and substitute $u_{k-1}$ by $u_0$. Then we obtain the linear
  approximate equation
  \begin{equation}
    \left(-\partial_t^2 + \partial_r^2 +\frac{2}r \partial_r + 5u_0^4\right) \eps  +
    e_{k-1} \approx 0
    \label{eps}  \end{equation}
  For $r \ll t$ we expect the time
  derivative to play a lesser role  so we neglect it and we are left
  with an elliptic equation with respect to the variable $r$,
  \begin{equation}
    \left(\partial_r^2 +\frac{2}r \partial_r + 5u_0^4\right) \eps  + e_{k-1} \approx 0, \qquad r
    \ll t
    \label{epsodd}
    \end{equation}
    For $r \approx t$ we can
  approximate $u_0^4$ by zero and rewrite \eqref{eps} in the form
 \begin{equation}
  \left(-\partial_t^2 + \partial_r^2 +\frac{2}r \partial_r \right) \eps +
  e_{k-1}\approx 0
  \label{epseven}\end{equation}
  Here the time and spatial derivatives have the same strength.
  However, we can identify another principal variable, namely $a =
  r/t$ and think of $\eps$ as a function of $(t,a)$.  Later, we reduce the above equation to a
  Sturm-Liouville problem in~$a$ which becomes singular at~$a=1$.

  The above heuristics lead us to a two step iterative construction of
  the $v_k$'s. The two steps successively improve the error in the two
  regions $r \ll t$, respectively $r \approx t$. To be precise, we
  define $v_k$ by
  \begin{equation}
    \left(\partial_r^2 +\frac{2}r \partial_r + 5u_0^4\right) v_{2k+1} +
    e_{2k}^0=0
    \label{vkodd}\end{equation} respectively
  \begin{equation}
    \left(-\partial_t^2+\partial_r^2 +\frac{2}r \partial_r \right) v_{2k} +
    e_{2k-1}^0=0
    \label{vkeven}\end{equation} both equations having zero Cauchy
  data\footnote{The coefficients are singular at $r=0$, therefore this
    has to be given a suitable interpretation} at $r=0$.  Here at each
  stage the error term $e_k$ is split into a principal part and a
  higher order term (to be made precise below),
  \[
  e_k = e_k^0 + e_k^1
  \]
  The successive errors are then computed as
  \[
  e_{2k} = e_{2k-1}^1 + N_{2k} (v_{2k}), \qquad e_{2k+1} = e_{2k}^1 -
  \partial_t^2 v_{2k+1} + N_{2k+1} (v_{2k+1})
  \]
  where
  \begin{equation}
    N_{2k+1}(v) =5(u_{2k}^4-u_0^4) \,v +
  10 u_{2k}^3 \,v^2 + 10 u_{2k}^2 \,v^3 + 5u_{2k} \,v^4 + v^5
    \label{eodd}\end{equation}
    respectively
  \begin{equation}
    N_{2k}(v) = 5 u_{2k-1}^4v +
  10 u_{2k-1}^3 v^2 + 10 u_{2k-1}^2 v^3 + 5u_{2k-1} v^4 + v^5
    \label{eeven}\end{equation}
  To formalize this scheme we need to introduce suitable function
  spaces in the cone \[\calC_0=\{(r,t)\::\:0\le r<t,\, 0<t<t_0\}\] for
  the successive corrections and errors. We first consider the $a$
  dependence. For the corrections $v_k$ we
set
\[
\drnu = \frac{\nu-1}2>-\frac12
\]
and use

\begin{defi}\label{def:Q} For $i \in \Nat$ we let $j(i)=0$ if $\nu$ is
  irrational, respectively $j(i) = i$ if $\nu$ is rational.

a)  For any positive integer $k$, we define $\mathcal Q$ to be the algebra of continuous functions $q:[0,1] \to
  \R$ with the following properties:

  (i) $q$ is analytic in $[0,1)$ with an even expansion at $0$.

  (ii) Near $a=1$ we have an absolutely convergent expansion of the
  form
  \[
\begin{split}
  q(a) = q_0(a) + \sum_{i =1}^\infty  (1-a)^{i (\drnu+1) - 2\left[\frac{i-1}4\right]}\sum_{ j=0}^{ j(i)}
    q_{ij}(a) (\ln (1-a))^j
  \end{split}
  \]
  with analytic coefficients $q_0$, $q_{ij}$.

  b) $\mathcal Q_m$ is the algebra which is defined similarly, with
  the additional requirement that
  \[
  q_{ij}(1) =0 \text{\ \ if\ }\ \ i =4k+1 \geq 4m+1.
  \]
\end{defi}

We remark that the exponents of $1-a$ in the above series are all
positive because of $\beta_0>-\frac12$.  For the errors $e_k$ we
introduce

\begin{defi}
   $\mathcal Q'$ is the space of continuous functions $q:[0,1) \to
  \R$ with the following properties:

  (i) $q$ is analytic in $[0,1)$ with an even expansion at $0$

  (ii) Near $a=1$ we have a convergent expansion of the form
  \[\begin{split}
  q(a) = q_0(a) + \sum_{i =1}^\infty  (1-a)^{i (\drnu+1) - 2\left[\frac{i-1}4\right]-1}\sum_{ j=0}^{ j(i)}
    q_{ij}(a) (\ln (1-a))^j
  \end{split}
  \]
  with analytic coefficients $q_0$, $q_{ij}$.

  b) $\mathcal Q'_m$ is the space which is defined similarly, with
  the additional requirement that
  \[
 q_{ij}(1) =0 \text{\ \ if\ }\ \ i =4k+1 \geq 4m+1.
  \]
\label{def:Qp}\end{defi}

By construction, $\mathcal Q_k\subset \mathcal Q'_k$. The families
$\cQ'$ and $\cQ_k'$ are obtained by applying $a^{-1}\partial_a$ to
the algebras $\cQ$ and $\cQ_k$, respectively.

We remark that the number of logarithms in these definitions in the
case when $\nu$ is rational is far from optimal, but we have chosen
this form since it simplifies the presentation. Next we define the
class of functions of $R$:

\begin{defi}
  $S^m(R^k)$ is the class of analytic functions
  $v:[0,\infty) \to \R$ with the following properties:

  (i) $v$ vanishes of order $m$  and $R^{-m} v$ has an even
  Taylor expansion at $R=0$.

  (ii) $v$ has a convergent expansion near $R=\infty$,
  \[
  v = \sum_{i=0}^\infty c_{i}\, R^{k-2i}
  \]
\end{defi}

The importance of even expansions in $R$ lies with the fact that
only those correspond to smooth functions in~$\R^3$. For the same
reason, we will work with even~$m$.  We also introduce another
auxiliary variable,
\begin{equation}\label{eq:bdef} b = \frac{1}{(t \lambda)^2}
\end{equation}
Since we seek solutions inside the cone we can restrict $b$ to a
small interval $[0,b_0]$. We combine these three components in order
to obtain the full function class which we need:

\begin{defi}
  a) $S^m(R^k ,\mathcal Q_n)$ is the class of analytic
  functions $v:[0,\infty) \times [0,1]\times [0,b_0] \to \R$ so that

  (i) $v$ is analytic as a function of $R,b$,
  \[
  v: [0,\infty) \times [0,b_0] \to \mathcal Q_n
  \]

  (ii) $v$ vanishes of order $m$ and $R^{-m} v$ has an even
  Taylor expansion at $R=0$.

  (iii) $v$ has a convergent expansion at $R=\infty$,
  \[
  v(R,\cdot,b) = \sum_{i=0}^\infty c_{i}(\cdot,b) R^{k-2i}
  \]
  where the coefficients $c_{i}: [0,b_0] \to \cQ_m$ are analytic with
  respect to $b$.

  b) $\IS^m(R^k ,\mathcal Q_n)$ is the class of analytic
  functions $w$ on the cone $\calC_0$ which can be represented as
  \[
  w(r,t) = v(R,a,b), \qquad v \in S^m(R^k ,\mathcal Q_n)
  \]
\end{defi}

We note that the representation of functions on the cone as in part
(b) is in general not unique since $R,a,b$ are dependent variables.
Later we shall exploit this fact and switch from one representation
to another as needed. We shall prove by induction that the
successive corrections $v_k$ and the corresponding error terms $e_k$
can be chosen with the following properties: For each $k\ge1$,

\begin{align}
  v_{2k-1} &\in \frac{\lambda^{\frac12}}{(t \lambda)^{2k}} \IS^2(R,\cQ_{k-1})
  \label{v2k-1}\\
  t^2 e_{2k-1} &\in \frac{\lambda^{\frac12}}{(t \lambda)^{2k}} \IS^0(R, \cQ'_{k-1})
  \label{e2k-1}\\
  v_{2k} &\in \frac{\lambda^{\frac12}}{(t \lambda)^{2k+2}} \IS^2(R^3,\cQ_{k}) \label{v2k}\\
  t^2 e_{2k} &\in \frac{\lambda^{\frac12}}{(t \lambda)^{2k}} \big[\IS^0(R^{-1} ,\cQ_k)
   + b\IS^0(R ,\cQ'_{k})    \big]
  \label{e2k}\end{align}

\noindent with \eqref{e2k} also valid for $k=0$.  We remark that the
order of vanishing at $R=0$ can be successively improved with $k$,
but this does not appear to be important.
\medskip

{\bf Step 0:} {\em The analysis at $k=0$}

\medskip
\noindent With our notations, one checks that
\begin{equation}\label{e20}
t^2 e_0 = -t^2\partial_{tt}[ \lambda^{\frac12} W(\lambda(t) r)] \in
\lambda^{\frac12} IS^0(R^{-1})
\end{equation}
as claimed. Now assume we know the above relations hold up to $k-1$
with $k\ge1$, and we show how to construct $v_{2k-1}$, respectively
$v_{2k}$, so that they hold for the index $k$.

\medskip

{\bf Step 1:} {\em Begin with $e_{2k-2}$ satisfying \eqref{e2k}
  or~\eqref{e20} and choose $v_{2k-1}$ so that \eqref{v2k-1} holds.}

\medskip \noindent If $k=1$, then define $e_0^0:= e_0$. If $k>1$,
we use \eqref{e2k} to write
\[
  e_{2k-2} = e_{2k-2}^0 + e_{2k-2}^1
\]
where
\[
t^2 e_{2k-2}^0 \in \frac{\lambda^{\frac12}}{(t \lambda)^{2k-2}}
\IS^0(R^{-1} ,\cQ_{k-1}), \qquad t^2 e_{2k-2}^1 \in
\frac{\lambda^{\frac12}}{(t \lambda)^{2k}} \IS^0(R ,\cQ'_{k-1})
\]
In the first term we can set $b=0$ and eliminate the $b$ dependence,
as all the $b$ dependent part can be included in the second term.

We note that the term $e_{2k-2}^1$ can be included in $e_{2k-1}$,
cf.~\eqref{e2k-1}. We define $v_{2k-1}$ as in \eqref{vkodd}
neglecting the $a$ dependence of $e_{2k-2}^0$. In other words, $a$
is treated as a parameter.  Changing variables to $R$ in
\eqref{vkodd} we need to solve the equation
% \footnote{When $k>1$ the right-hand side here can be
%   improved to $t^2 e_{2k-2}^0 \in \frac{\lambda^{\frac12}}{(t
%     \lambda)^{2k}} \IS^1(R^{-1}, \cQ_{k-1})$, but not for $k=1$.}
\[
(t \lambda)^2 L v_{2k-1} = t^2 e_{2k-2}^0 \in
\frac{\lambda^{\frac12}}{(t
  \lambda)^{2k-2}} \IS^0(R^{-1}, \cQ_{k-1} )
\]
where the operator $L$ is given by
\[
L = -\partial_R^2 -\frac{2}R \partial_R - 5 W^4(R)
\]
Then \eqref{v2k-1} is a consequence of the following ODE lemma.

\begin{lemma} The solution $v$ to the equation
  \[
  L v= f \in S^0(R^{-1}), \qquad v(0) = v'(0)=0
  \]
  has the regularity
  \[
  v \in S^2(R)
  \]
\end{lemma}

\begin{proof}
  Since $f$ is analytic at $0$ with a constant leading term, one can
  easily write down an even Taylor series for $v$ at $0$ with a quadratic
  leading term.

  It remains to determine the asymptotic behavior of $v$ at infinity.
  For this it is convenient to remove the first order derivative in
  $L$ (to achieve constancy of the Wronskian). Thus, we seek a
  solution of
  \[ \tilde L {R}\, v= {R} f, \qquad \tilde L= \partial_R^2 + 5W^4=\partial_R^2 + \frac{5}{(1+R^2/3)^2}
  \]
  We use this fundamental system of solutions for $\tilde L$:
  \[\begin{split}
  \phi(R) &= R(1-R^2/3)(1+R^2/3)^{-\frac32}\\ \theta(R) &=
  (1+R^2/3)^{-\frac32} (1-2R^2+R^4/9)
  \end{split}
  \]
  Clearly, $L\partial_\lambda W=0$ and we set $\phi= R \partial_\lambda
  W|_{\lambda=1}$ up to a constant. The function $\theta$ is then
  determined from the  Wronskian constancy condition $W(\theta,\phi)=1$.
  This allows us to obtain an
  integral representation for $v$ using the variation of parameters
  formula, which gives
  \[
  v = -R^{-1}\theta(R)\int_0^R \phi(R')R' f(R')\,
  \,dR' +  R^{-1} \phi(R) \int_0^R \theta(R') {R'}
  f(R')\, \,dR'
  \]
The right-hand side grows
  like $R$, as claimed.
\end{proof}

\noindent As a special case of the above computation we note the
representation for $v_1$,
\begin{equation}\label{eq:v1}
  v_1 = \frac{\lambda^{\frac12}}{(t \lambda)^{2}} V(R), \qquad V \in S^2(R)
\end{equation}
This justifies the choice of the second term in the expansion for
$u_{2k-1}$ in Theorem~\ref{thm:sec2}.

\medskip

{\bf Step 2:} {\em Show that if $ v_{2k-1}$ is chosen as above then
  \eqref{e2k-1} holds.}

\medskip \noindent Thinking of $v_{2k-1}$ as a function of $t$, $R$
and $a$ we can write $e_{2k-1}$ in the form
\[
e_{2k-1} = N_{2k-1}(v_{2k-1}) + E^t v_{2k-1} + E^a v_{2k-1}
\]
Here $N_{2k-1}(v_{2k-1})$ accounts for the contribution from the
nonlinearity and is given by \eqref{eodd}. $E^t v_{2k-1}$ contains
the terms in
\begin{equation}\label{eq:prtt}
  -\partial_{tt} v_{2k-1} (t,R,a)
\end{equation}
where no derivative applies to the variable $a$, while $ E^a
v_{2k-1}$ contains those terms in
\[
\Big(\partial_{tt}-\pr_{rr} -\frac{2}{r}\pr_r\Big) v_{2k-1}(t,R,a)
\]
where at least one derivative applies to the variable~$a$ (recall
that in Step~1 the parameter $a$ was frozen). We begin with the
terms in $N_{2k-1}$. We first note that, by summing the $v_j$ over
$1\le j\le 2k-2$,
\begin{equation}\label{eq:uentw1}
  u_{2k-2} - u_0 \in \frac{\lambda^{\frac12}}{(t \lambda)^{2}} \IS^2(R,\mathcal
  Q_{k-1})\end{equation}
The first term in $N_{2k-1}(v_{2k-1})$ contributes
\begin{eqnarray}
t^2 (u_{2k-2}^4-u_0^4)v_{2k-1}  &=& t^2 [(u_{2k-2}-u_0)^4 +
4(u_{2k-2}-u_0)^3u_0 \nn\\
&&\qquad + 6(u_{2k-2}-u_0)^2 u_0^2 + 4 (u_{2k-2}-u_0) u_0^3]
v_{2k-1}\label{eq:udiff1}
\end{eqnarray}
Using \eqref{eq:uentw1} we compute
\begin{align*}
t^2 (u_{2k-2}-u_0)^4 v_{2k-1} &\in \frac{1}{(t\lambda)^6}
\IS^8(R^4,\cQ_{k-1}) \frac{\lambda^{\frac12}}{(t \lambda)^{2k}}
\IS^2(R,\cQ_{k-1}) \\
&\subset a^6 \IS^2(R^{-2},\cQ_{k-1}) \frac{\lambda^{\frac12}}{(t
\lambda)^{2k}} \IS^2(R,\cQ_{k-1})\\
& \subset \frac{\lambda^{\frac12}}{(t \lambda)^{2k}}
\IS^2(R^{-1},\cQ_{k-1})
\end{align*}
as well as
\begin{align*}
t^2 (u_{2k-2}-u_0)u_0^3 v_{2k-1} &\in t^2
\frac{\lambda^{\frac12}}{(t \lambda)^{2}} \IS^2(R,\mathcal
  Q_{k-1}) \lambda^{\frac32} S^0(R^{-3})
 \frac{\lambda^{\frac12}}{(t \lambda)^{2k}}
\IS^2(R,\cQ_{k-1})\\
&\subset \frac{\lambda^{\frac12}}{(t \lambda)^{2k}}
\IS^2(R^{-1},\cQ_{k-1})
\end{align*}
The other two terms in \eqref{eq:udiff1} are similar. Next, compute
\begin{align*}
  t^2 v_{2k-1}^5 &\in \frac{t^2 \lambda^{\frac52}}{(t
\lambda)^{10k}} \IS^{10}(R^5,\cQ_{k-1}) \\
&\subset \frac{\lambda^{\frac12}R^6}{(t \lambda)^{10k-2}}
\IS^4(R^{-1},\cQ_{k-1})\\
&\subset \frac{\lambda^{\frac12}}{(t \lambda)^{2k}} {a^6} b^{4(k-1)}
\IS^2(R^{-1},\cQ_{k-1}) \subset \frac{\lambda^{\frac12}}{(t
\lambda)^{2k}}
 \IS^2(R^{-1},\cQ_{k-1})
\end{align*}
and
\begin{align*}
  t^2 u_{2k-2}^3\,v_{2k-1}^2 &\in
\lambda^{-\frac12} (t\lambda)^2 \IS^0(R^{-3}, \cQ_{k-1})
\frac{\lambda}{(t\lambda)^{4k}} \IS^4(R^2,\cQ_{k-1}) \\
&\subset \frac{\lambda^{\frac12}}{(t\lambda)^{2k}} b^{2k-2}
\IS^4(R^{-1},\cQ_{k-1})\subset
\frac{\lambda^{\frac12}}{(t\lambda)^{2k}}  \IS^2(R^{-1},\cQ_{k-1})
\end{align*}
with similar statements for $u_{2k-2}^2 v_{2k-1}^3$ and $u_{2k-2}
v_{2k-1}^4$.
 Summing up we obtain
\[
N_{2k-1}(v_{2k-1}) \in \frac{\lambda^{\frac12}}{(t\lambda)^{2k}}
\IS^2(R^{-1},\cQ_{k-1}) \subset
\frac{\lambda^{\frac12}}{(t\lambda)^{2k}} \IS^2(R^{-1},\cQ_{k-1}')
\]
This concludes the analysis of $N_{2k-1}(v_{2k-1})$. We continue
with the terms in $E^t v_{2k-1}$, where we can neglect the $a$
dependence. Therefore, it suffices to compute
\[
t^2 \partial_t^2 \left( \frac{\lambda^{\frac12}}{(t \lambda)^{2k}}
\IS^2(R)\right) \subset \frac{\lambda^{\frac12}}{(t \lambda)^{2k}}
\IS^2(R)
\]
Finally, we consider the terms in $E^a v_{2k-1}$. With
\[
v_{2k-1}(r,t) = \frac{\lambda^{\frac12}}{(t \lambda)^{2k}} w(R,a),
\quad w \in S^2(R,\mathcal Q_{k-1})
\]
we have
\begin{eqnarray*}
  t^2 E^a v_{2k-1}& =&   -2t\partial_t\left(\frac{\lambda^{\frac12}}{(t \lambda)^{2k}}\right)
  aw_a(R,a) + \frac{\lambda^{\frac12}}{(t \lambda)^{2k}}
  \big[ 2(\nu+1)aR w_{aR}(R,a)   \\
  && - 2Ra^{-1} w_{Ra}  - 2a^{-1} w_a(R,a)+
  (a^2-1) w_{aa}(R,a) + 2a w_a(R,a)\big]
\end{eqnarray*}
Since $\mathcal Q_{k-1}$ are even in $a$ we conclude that
\[
(1-a^2) \partial_{aa},\, a \partial_a,\, a^{-1}\pr_a\,:\; \mathcal
Q_{k-1} \to \mathcal Q'_{k-1}
\]
and therefore
\[
t^2 E^a v_{2k-1} \in \frac{\lambda^{\frac12}}{(t \lambda)^{2k}}
\IS^2(R,\mathcal Q'_{k-1})
\]
This concludes the proof of \eqref{e2k-1}. We remark that for the
special case of $k=1$, i.e., with $v_1$ as in~\eqref{eq:v1}, these
arguments yield
\begin{equation}
  \label{eq:e1_def}t^2 e_1 \in \frac{\lambda^{\frac12}}{(t\lambda)^{2}} IS^2(R)
\end{equation}

\medskip

{\bf Step 3:} {\em Define $v_{2k}$ so that \eqref{v2k} holds.}

\medskip\noindent We begin the analysis with $e_{2k-1}$ replaced by its main
asymptotic component $\tilde e^0_{2k-1}$ around $R = \infty$.
This has the form
\begin{equation}
  t^2 \tilde e^0_{2k-1} = \frac{\lambda^{\frac12}R}{(t \lambda)^{2k}}  q(a),
\qquad q \in \mathcal Q'_{k-1} \label{eq:e2kmin1}
\end{equation}
which we rewrite as
\[
t^2 \tilde e^0_{2k-1} = \frac{\lambda^{\frac12}}{(t \lambda)^{2k-1}}  a
q(a)
\]
We remark that \eqref{eq:e1_def} implies that \[t^2 \tilde e^0_1(a)=
\frac{\lambda^{\frac12}}{t\lambda} a.\]  Consider the equation
\eqref{vkeven} with $ \tilde e^0_{2k-1}$ on the right-hand side,
\[
t^2\left(-\partial_t^2+\partial_r^2 +\frac{2}r \partial_r\right)
\tilde v_{2k} = -t^2 \tilde e^0_{2k-1}
\]
We look for a solution $\tilde v_{2k}$ which has the form
\[
\tilde v_{2k} = -\frac{\lambda^{\frac12}}{(t \lambda)^{2k-1}}
W_{2k}(a)
\]
Thus,
\[
\begin{split}
  t^2 \left(-\partial_t^2+\partial_r^2 +\frac{2}r \partial_r \right) \left(\frac{\lambda^{\frac12}}{(t \lambda)^{2k-1}}
    W_{2k}(a)\right) = \frac{\lambda^{\frac12}}{(t \lambda)^{2k-1}} a q(a)
\end{split}
\]
Conjugating out the power of $t$ we get
\[
t^2 \Big( -\Big(\partial_t+\frac{(2k-\frac32)\nu-\frac12}t\Big)^2+
  \partial_r^2 + \frac{2}r \partial_r
\Big)W_{2k}(a) = a q(a)
\]
which we rewrite as an equation in the $a$ variable,
\begin{equation}
  L_{(2k-\frac32)\nu-\frac12}\, W_{2k}(a) = a q(a), \qquad q \in \cQ'_{k-1}
  \label{w2kj}
  \end{equation}
  where the one-parameter family of
operators $L_{\beta}$ is defined by
\begin{equation}\label{eq:Lbeta_def}
L_\beta =(1-a^2) \partial_{aa} + 2(a^{-1} + a \beta - a)
\partial_a -\beta^2 + \beta
\end{equation}
We claim that solving this equation with zero Cauchy data at $a=0$
yields a solution which satisfies
\begin{equation}
  W_{2k}(a) =  a^3 q_1(a),        \qquad q_1  \in        \cQ_k
  \label{cqk}\end{equation}
This gives
\[
\tilde v_{2k} = \frac{\lambda^{\frac12}}{(t \lambda)^{2k-1}} a^3
q_1(a) = \frac{\lambda^{\frac12}}{(t \lambda)^{2k+2}} R^3
q_1(a)
\]
which is not entirely suitable as the next correction since it has an
odd expansion at $R=0$ instead of an even one. To remedy this
we simply set
\[
 v_{2k} = \frac{\lambda^{\frac12}}{(t \lambda)^{2k+2}} R^3 \frac{R}{(1+R^2)^\frac12}
q_1(a) \in \frac{\lambda^{\frac12}}{(t \lambda)^{2k+2}} S^4(R^3,\cQ_k)
\]
Clearly, this will conclude Step~3.   To prove the claim \eqref{cqk} we need the
following

\begin{lemma}
  Let $f \in a \cQ'_{k-1}$, $k\ge1$.
  Then there is a unique solution $w \in a^3 \cQ_k$ to the equation
  \begin{equation}
    L_{(2k-\frac32)\nu-\frac12}\, w = f, \qquad w(0) = 0, \ \partial_a w(0) = 0
  \end{equation}
\end{lemma}

\begin{proof}
  Denote
  \begin{equation}\label{eq:b_def} \beta= \big(2k-\frac32
  \big)\nu-\frac12 = (4k-3)\beta_0+2(k-1)>-\frac12.
  \end{equation}
  % Since $k\ge1$ and $\nu>1$, also
  %$\beta>0$.
  We write
\[
L_\beta =a^{-2} \partial_a (a^2 \partial_a)-a^2 \partial_{aa} + 2(
\beta - 1)a\partial_a -\beta^2 + \beta
\]
To study the behavior of the solutions at $0$  we match the coefficients in
  $L_\beta \, w=f$ with
  \[
  f(a) = \sum_{j=1}^\infty f_j\, a^{2j-1}, \qquad w(a) = \sum_{j=2}^\infty w_j \, a^{2j-1}
  \]
  yields the system, with $j\ge1$,
  \[
  (2j+1)(2j+2) w_{j+1} + [-(2j-1)(2j-2) + 2(\beta-1) (2j-1) +
  (\beta-\beta^2)]w_j=f_j
  \]
  where we take $w_1=0$. The coefficient of $w_{j+1}$ is always nonzero;
  this allows us to successively compute the coefficients $w_j$. The
  convergence of the series for $w$ follows from the
  convergence of the series for $f$.

   It remains to study the solution $w$ near $a=1$. The behavior of
  $L_{\beta}$ at $1$ is well approximated by
  \[
  L_\beta^1 = 2(1-a) \partial_{aa} + 2\beta  \partial_a =
  2(1-a)^{\beta+1}\partial_a  [(1-a)^{-\beta}\partial_a]
\]
which annihilates the functions $1$ and $(1-a)^{\drb+1}$. Therefore,
we seek a fundamental system for $L_\beta  y=0$ of the form
  \begin{equation}\label{eq:fundsys}
  \phi_1(a) = 1+\sum_{\ell=1}^\infty \mu_\ell (1-a)^\ell , \qquad \phi_2(a) =
  (1-a)^{\drb+1}\Big[ 1+\sum_{\ell=1}^\infty \tilde \mu_\ell (1-a)^\ell \Big]
  \end{equation}
This leads to the conditions, with $\mu_0=\tilde\mu_0=1$,
\begin{align}
\label{eq:phi1def} 2\mu_{\ell+1}(\ell+1)(\ell-\beta) &= \mu_\ell
[\ell(\ell-2\beta+3)+\beta^2-\beta] + 2\sum_{j=0}^{\ell-1} j\mu_j \\
2\tilde\mu_{\ell+1}(\ell+\beta+2)(\ell+1) &= \tilde \mu_\ell
[(\beta+1+\ell)(\ell-\beta+4)+\beta^2-\beta ] \label{eq:phi2def}\\
&\qquad + 2\sum_{j=0}^{\ell-1} (\beta+1+j)\tilde\mu_j
 \nonumber
\end{align}
Clearly, \eqref{eq:phi2def} always has a solution whereas
\eqref{eq:phi1def} requires $\beta\not\in \Z_0^+$; in the latter
case, the series in~\eqref{eq:fundsys} define entire functions. If,
on the other hand, $\beta\in \Z_0^+$, then
  $\phi_1$ is modified to
  \begin{equation}\label{eq:phi1mod}
  \phi_1(a) = 1+\sum_{\ell=1}^\infty \mu_\ell (1-a)^\ell + c_1\phi_2(a) \ln(1-a)
  \end{equation}
with some unique choice of $c_1$.

Modulo a linear combination of $\phi_1, \phi_2$ it suffices to find
one solution to the inhomogeneous equation $L_\beta\, w=f$ near
$a=1$. We begin with the case when $f$ has the form
\[
f(a) = (1-a)^\gamma  \sum_{k=0}^\infty f_k (1-a)^k
\]
with $\gamma > 0$.  Then we seek a solution $w$ of the form
\[
w(a)= (1-a)^{\gamma+1}  \sum_{k=0}^\infty w_k (1-a)^k
\]
This leads to the system
\begin{equation}
\left\{
  \begin{array}{ll}
    2 w_0 (\gamma+1)(\gamma -\beta) &= f_0 \\
    \displaystyle 2 w_{\ell+1}
    (\gamma+\ell+2)(\gamma-\beta+\ell+1) &=  f_\ell + w_\ell[(\gamma+\ell+1)(\gamma+\ell+4-2\beta)+\beta^2
    -\beta] \\
     &\quad +2 \sum_{j=0}^{\ell-1} (\gamma+j+1)  w_j
\end{array}\right.
\label{eq:phi2defb}
\end{equation}
which is solvable unless $\beta - \gamma$ is a nonnegative integer,
in which case the solution $w$ has the modified form
\[
w(a)= (1-a)^{\gamma+1}  \sum_{k=0}^\infty w_k (1-a)^k + c \phi_2(a) \ln(1-a)
\]
Indeed, assume that $\beta-\gamma=p\ge 0$ with integer~$p$. Then
\eqref{eq:phi2defb} is in general violated for $\ell+1=p$ regardless
of the choice of $w_p$; thus, $c$ above has to be chosen so that the
coefficients of $(1-a)^{\gamma+p}=(1-a)^\beta$ in $L_\beta w=f$
match. This can be done since
\[
\phi_2(a)\ln(1-a) = [(1-a)^{\beta+1} + O((1-a)^{\beta+2})] \ln(1-a)
\]
and
\[ L_\beta [\phi_2(a)\ln(1-a) ] = 2(\beta+1) (1-a)^{\beta} +
O((1-a)^{\beta+1})
\]
with $\beta+1>\frac12$ and thus nonzero. Note that $w_p$ remains
undetermined -- this amounts to the freedom of adding a constant
multiple of $\phi_2(a)$ to $w(a)$.

Similarly, if $f$ has the form
\[
f(a) = \sum_{m=0}^j (\ln(1-a))^m f_m(a) =\sum_{m=0}^j (\ln(1-a))^m
(1-a)^\gamma \sum_{k=0}^\infty f_{km} (1-a)^k
\]
then we seek a solution $w$ of the form
\[
w(a)= \sum_{m=0}^j (\ln(1-a))^m w_m(a) = \sum_{m=0}^j (\ln(1-a))^m
(1-a)^{\gamma+1} \sum_{k=0}^\infty w_{km} (1-a)^k
\]
Identifying the coefficients of the powers of $\ln(1-a)$ we obtain
the system
\[
L_\beta w_j = f_j, \qquad L_\beta w_{j-1} = f_{j-1}- j Q_1 w_j,
\]
respectively
\[
L_{\beta} w_{j-k} = f_{j-1}- j Q_1 w_{j-k+1} + j(j-1) Q_2 w_{j-k+2},
\qquad k \geq 2
\]
where
\[
Q_1 = 2(1+a)
\partial_a + \frac{1+a}{1-a} + \frac{2(1+a)}a + \frac{2a\beta}{1-a},
\qquad Q_2 = \frac{1+a}{1-a}
\]
This system is solved iteratively as in the first case provided that
$\beta-\gamma$ is not a nonnegative integer. Otherwise,
the solution $w$ has the modified form
\[
w(a)=    \sum_{m=0}^j (\ln(1-a))^m
\big[c_m\,\phi_2(a)\ln(1-a)+(1-a)^{\gamma+1} \sum_{k=0}^\infty
w_{km} (1-a)^k \big]
\]

Consider now $f \in \cQ'_{k-1}$.  If $\nu$ is irrational then
according to Definition~\ref{def:Qp}
we can represent it in the form
\[
f(a) = f_0(a) + \sum_{i =1}^\infty (1-a)^{i (\drnu+1) -
  2\left[\frac{i-1}4\right]-1} f_{i}(a)
\]
with $f_{i}$ analytic. Hence the exponent $\gamma$ above takes
the values
\[
\gamma_i = i (\drnu+1) - 2\left[\frac{i-1}4\right]-1
\]
On the other hand, we have
\[
\beta = (2k - \frac32)\nu -\frac12 = (4k-3) (\drnu+1) - 2k+1 =
\gamma_{4k-3}
\]
Then $\gamma_i-\beta$ can only be an integer if $i=4k-3$. However,
in this case there is no logarithmic term due to the additional
condition $q_{4k-3}(1)=0$ which has the effect of replacing $
\gamma_{4k-3}$ by $\gamma_{4k-3}+1$. The conclusion of the lemma
follows in the irrational case.

On the other hand, if $\nu$ is rational then $f$ has the
representation
\[
f(a) = f_0(a) + \sum_{i =1}^\infty (1-a)^{i (\drnu+1) -
  2\left[\frac{i-1}4\right]-1}\sum_{ j=0}^{i} f_{ij}(a) (\ln
(1-a))^j
\]
with $f_{ij}$ analytic.  In order for $\beta-\gamma_\ell$ to be a
nonnegative integer we need to have $\ell \leq 4k-3$.  Again the
condition $q_{4k-3}(1)=0$ guarantees that there is no extra
logarithmic term arrising from the $i = 4k-3$ component of $f$. On
the other hand, if $\ell < 4k-3$ then we can contribute an extra
logarithm for a total of $\ell+1 \leq  4k-3$ logarithms to the $i =
4k-3$ term of $w$. The conclusion of the lemma again follows.

It is also worth  discussing the special case $k=1$. This will also
serve to explain how the algebra $\cQ_k$ arises  in the iteration.
If $k=1$,  then~\eqref{w2kj} reduces to the equation,  with the
usual $\beta_0=(\nu-1)/2$,
\[
L_{\beta_0} \, W_2(a)=a\] due to $t^2 f_1(a)=
\lambda^{\frac12}(t\lambda)^{-1} a$. As discussed above,
\[\begin{split}
W_2(a) &= g_0(a) +  g_1(a)(1-a)^{\drnu+1} \text{\ \ if\ }\beta_0\not\in\Z_0^+\\
 W_2(a) &= h_0(a)
+ h_1(a)(1-a)^{\drnu+1} + h_2(a)(1-a)^{\drnu+1}\log(1-a) \text{\ \
if\ }\beta_0\in\Z_0^+
\end{split}
\]
Thus, we see that in all cases $W_2\in \cQ_1$ for $j=0,1$ and $a$
near~$1$.
\end{proof}

\medskip
{\bf Step 4}: {\em With $v_{2k}$ as above show that $e_{2k}$ is as
claimed.}

\medskip
\noindent
Modifying \eqref{eq:e2kmin1} to insure an even expansion at $R=0$ we set
\[
t^2 e_{2k-1}^0 = \frac{\lambda^\frac12}{(t\lambda)^{2k}}
R \frac{R}{(R^2+1)^\frac12} q(a), \qquad q \in \cQ_{k-1}'
\]
Then we can write $e_{2k}$ in the form
\[\begin{split}
t^2 e_{2k} &= t^2 \left(e_{2k-1} - e_{2k-1}^0 \right) + t^2 \left(
  e_{2k-1}^0- \left(-\partial_t^2+\partial_r^2 + \frac{2}r \partial_r\right) v_{2k} \right) + t^2
  N_{2k}(v_{2k})
\end{split}
\]
where $N_{2k}$ is defined by \eqref{eeven}.

We begin with the first term in $e_{2k}$, which has the form
\[
t^2(e_{2k-1} - e_{2k-1}^0) \in \frac{\lambda^{\frac12}}{(t
\lambda)^{2k}} \IS^0(R^{-1},\mathcal Q'_{k-1})
\]
 We claim that
\begin{equation}
  \IS^0(R^{-1},\mathcal Q'_{k-1}) \subset \IS^0(R^{-1}
  ,\mathcal Q_{k-1}) + \frac{1}{(t\lambda)^{2}}  \IS^2(R,\mathcal
  Q'_{k-1}) \label{qqp}
  \end{equation}
  For $w \in \IS^1(R^{-1},\mathcal Q'_{k-1})$ we write
\[
w = (1-a^2)w + \frac{1}{(t \lambda)^2} R^{2} w
\]
Then
\[
(1-a^2)w \in \IS^0(R^{-1},\mathcal Q_{k-1}), \qquad \frac{1}{(t
\lambda)^2} R^{2} w \in   \frac{1}{(t\lambda)^{2}}  \IS^2(R,\mathcal Q'_{k-1})
\]
as desired.

Next, we consider the expression
\[
f = t^2 \left(e_{2k-1}^0- \left(-\partial_t^2+\partial_r^2 + \frac{2}r
    \partial_r\right) v_{2k} \right)
\]
By construction this would vanish if the factor $R
(1+R^2)^{-\frac12}$ were dropped in both $e_{2k-1}^0$ and $v_{2k}$.
Hence $f$ contains only the components of the second term for which
at least one derivative falls on the factor $R (1+R^2)^{-\frac12}$:
\[
\begin{split}
f &= \tilde v_{2k} t^2 \left(-\partial_t^2+\partial_r^2 + \frac{2}r
    \partial_r\right) \frac{R}{(1+R^2)^\frac12}  - 2 t^2 \partial_t\tilde
  v_{2k}  \partial_t \frac{R}{(1+R^2)^\frac12} \\ & + 2 t^2 \partial_r\tilde
  v_{2k}  \partial_r \frac{R}{(1+R^2)^\frac12} + \tilde
  v_{2k} \frac{2 t^2} r \partial_r \frac{R}{(1+R^2)^\frac12}
\end{split}
\]
where $v_{2k}=\tilde v_{2k} R(1+R^2)^{-\frac12}$.  Then a direct
computation gives
\[
f \in \frac{\lambda^{\frac12}}{(t \lambda)^{2k+2}} IS^2(R,\cQ'_{k})
\]
as needed.

Finally we consider the nonlinear terms in
$N_{2k}(v_{2k})$. Again the $a,b$ dependence is uninteresting since
$\mathcal Q_k$ is an algebra. We start with the last term
in~\eqref{eeven}. By construction,
\[
v_{2k} \in \frac{\lambda^{\frac12}}{(t\lambda)^{2k}} \IS^2(R,\cQ_k)
\]
Thus,
\begin{align*}
  t^2 v_{2k}^5 &\in \frac{t^2\lambda^{\frac52}}{(t\lambda)^{10k}}
  \IS^5(R^5,\cQ_k) \subset
  \frac{\lambda^{\frac12}}{(t\lambda)^{2k}}\frac{1}{(t\lambda)^{8k-2}}(1+R^2)^3
  \IS^1(R^{-1},\cQ_k)
\end{align*}
Using that
\[
\frac{(1+R^2)^3}{(t\lambda)^{8k-2}} = b^{4k-1}+ {3a^2}b^{4k-2} +
3a^4 b^{4k-3} + a^6 b^{4k-4}
\]
we conclude that
\[
t^2 v_{2k}^5 \in \frac{\lambda^{\frac12}}{(t\lambda)^{2k}}
  \IS^1(R^{-1},\cQ_k)
\]
which is an admissible contribution to \eqref{e2k}. Finally, we
check the first term in \eqref{eeven}. From
\[
u_{2k-1}-u_0 \in \frac{\lambda^{\frac12}}{(t\lambda)^2}\IS^2 (R,
\cQ_k)
\]
and the form of $v_{2k}$,
\begin{align*}
t^2 u_{2k-1}^4 v_{2k} &\in t^2 \Big( \lambda^{\frac12} W(R) +
\frac{\lambda^{\frac12}}{(t\lambda)^2}\IS^2 (R, \cQ_k)\Big)^4
\frac{\lambda^{\frac12}}{(t\lambda)^{2k+2}}\IS^2 (R^3, \cQ_k)\\
&\subset t^2 \Big( \lambda^{\frac12} S^0(R^{-1}) + \lambda^{\frac12}
a^2 \IS^0 (R^{-1}, \cQ_k)\Big)^4
\frac{\lambda^{\frac12}}{(t\lambda)^{2k+2}}\IS^2 (R^3, \cQ_k)\\
& \subset (t\lambda)^2 \IS^0 (R^{-4}, \cQ_k)
\frac{\lambda^{\frac12}}{(t\lambda)^{2k+2}}\IS^2 (R^3, \cQ_k)\\
 &\subset \frac{\lambda^{\frac12}}{(t\lambda)^{2k}}\IS^2
(R^{-1}, \cQ_k)
\end{align*}
The other three terms in  \eqref{eeven} can be checked similarly.
 This concludes the proof of Theorem~\ref{thm:sec2}.
\end{proof}

\section{The linearized problem}
\label{sec:lin}

We seek a radial solution of \eqref{TheEquation} of the form
$u=u_{2k-1}+\eps$ with $u_{2k-1}$ as in Theorem~\ref{thm:sec2}. Our
ansatz leads to \beeq\label{eq:epsPDE}  \partial_{tt} \eps - \Delta
\eps - 5\lambda^2(t) W^4(\lambda(t)x) \eps = N_{2k-1}(\eps) +
e_{2k-1}
 \eneq where
$N_{2k-1}(\eps)$ is as in~\eqref{eodd} but with $u_{2k-2}$ replaced
with $u_{2k-1}$. Set $\eps(t,x)= v(\tau(t),\lambda(t) x)$ and
$y=\lambda(t) x$, $\tau=\tau(t)$. Then, with
$\dot\lambda=\f{d\lam}{d\tau}$,
\[ \partial_t \eps(t,r) = \tau'(t)(v_\tau + \dot\lambda \lambda^{-1} y\partial_y v) \]
and \beeq \label{eq:epsdiff} \partial_{tt} \eps(t,r) = \tau''(t)
(\partial_\tau + \dot\lambda \lambda^{-1} y\partial_y) v +
\tau'(t)^2 (\partial_\tau + \dot\lambda \lambda^{-1} y\partial_y)^2
v \eneq Set \[\tau(t) = \int_{t}^{t_0} \lambda(s)\,
ds+\frac{1}{\nu}t_0^{-\nu}= \frac{1}{\nu}t^{-\nu}\]
so that $\tau'(t)
= \lambda(t)$, and $\tau''(t)=\dot\lambda(\tau)\lambda(\tau)$ (we
are writing $\lambda(\tau)$ instead of $\lambda(t(\tau))$).
Then~\eqref{eq:epsPDE} can be rewritten as
\begin{equation}\label{eq:main1}\begin{split}
&[(\partial_\tau + \dot\lambda \lambda^{-1} y\partial_y)^2 v +
\dot\lambda \lambda^{-1} (\partial_\tau + \dot\lambda
\lambda^{-1} y\partial_y) v-\Delta v - 5W^4 v](\tau,y) \\
&= \lambda^{-2}(\tau)[N_{2k-1}(\eps) +
e_{2k-1}](t(\tau),\lambda^{-1}y)
\end{split}
\end{equation}
We remark that the linear Schr\"odinger operator $H=-\Delta -5W^4$
on $L^2(\R^3)$ has at least one negative eigenvalue as well as a
zero energy eigenvalue and resonance. The negative spectrum renders
the linear evolution in~\eqref{eq:main1} exponentially unstable. To
address this further, we switch to the radial variable $R=\lambda r$
from the previous section. Thus, \eqref{eq:main1} is the same as
\beeq\nn
\begin{split}
&[(\partial_\tau + \dot\lambda \lambda^{-1} R\partial_R)^2 v +
\dot\lambda \lambda^{-1} (\partial_\tau + \dot\lambda
\lambda^{-1} R\partial_R) v-v_{RR} -\frac{2}{R}v_R - 5W^4 v](\tau,R) \\
&= \lambda^{-2}(\tau)[N_{2k-1}(\eps) +
e_{2k-1}](t(\tau),\lambda^{-1}R)
\end{split}
\eneq or, with the new dependent variable $\tilde\eps(\tau,R)
:=Rv(\tau,R)$,
\beeq\label{eq:main2}
\begin{split}
&(\partial_\tau + \dot\lambda \lambda^{-1} R\partial_R)^2 \tilde\eps
- \dot\lambda \lambda^{-1} (\partial_\tau + \dot\lambda
\lambda^{-1} R\partial_R) \tilde\eps + \cL \tilde\eps \\
&= \lambda^{-2}R[N_{2k-1}(R^{-1}\tileps) + e_{2k-1}]
\end{split}
\eneq where
\[
\cL = -\pr_{RR} - 5W^4(R) \text{\ \ on\ \ } L^2(0,\infty)
\]
with a Dirichlet boundary condition at $R=0$.  Let
\[\beta(\tau):=\dot\lambda \lambda^{-1}(\tau)=\frac{1+\nu}{\tau\nu},
\quad \calD:=\pr_\tau + \beta(\tau)R\pr_R\] and
rewrite~\eqref{eq:main2} in the form
\begin{align}
& \calD^2\tileps-\beta(\tau)\calD \tileps + \cL \tileps = f
\label{eq:lin_eq}
\end{align}
To solve this equation we need precise spectral information
on the operator $\cL$.

\section{The spectral and scattering theory of the linearized
  operator}
\label{sec:spec}

\begin{defi}
  Let
  \[ \cL :=
  -\pr_{RR} - \frac{5}{(1+R^2/3)^2} \] be the half-line operator on
  $L^2(0,\infty)$ with a Dirichlet condition at $R=0$. It is self-adjoint
  on the domain
  \[ \dom(\cL)=\{ f\in L^2((0,\infty))\::\: f,f'\in
  AC([0,R])\,\forall \, R,\;f(0)=0,\, f'' \in L^2((0,\infty)) \}
  \]
\end{defi}

Note that $\cL\phi=0$ where \[\phi(R) :=
2R\partial_{\lambda}\Big|_{\lambda=1} \lambda^{\frac12} W(\lambda R)
= R(1-R^2/3)(1+R^2/3)^{-\frac32}
\]
This means that $\cL$ has a resonance at zero energy. Since $\phi$
has a single positive zero, it follows from oscillation theory, see
\cite{DS}, that there is an unique simple negative eigenvalue which
we denote by $\xi_d$. Thus, there is $\go\in L^2(0,\infty)\cap
C^\infty([0,\infty))$, decaying exponentially, and with $\go(R)>0$
for $R>0$ but $\go(0)=0$ so that $\cL \go = \etao \go$. We also
assume that $\|\go\|_2=1$. Clearly, $\cL$ has no other eigenvalues
or resonances.

\begin{lemma}
  \label{lem:spec} The spectrum of $\cL$ equals \[ \spec(\cL)= \{\etao\}\cup
  [0,\infty).\] The positive spectrum is purely
  absolutely continuous, and $\etao$ is a simple eigenvalue with
  eigenfunction that we denote by $\go$.
  Moreover, $\cL$ has a resonance at zero. In fact, $\cL\phi_0=0$ with
  $\phi_0(R)= R(1-R^2/3)(1+R^2/3)^{-\frac32}$. Finally, $\cL$ is in
  the limit-point case at infinity.
\end{lemma}

As usual, see Marchenko~\cite{Mar} or Section 2 of~\cite{GZ}, we
introduce the standard fundamental system of solutions $\phi(R,z)$
and $\theta(R,z)$ for all $z\in\Compl$ of $\cL y=zy$ with the
boundary conditions
\[ \phi(0,z)=\theta'(0,z)=0,\quad  \phi'(0,z)=\theta(0,z)=1 \]
so that in particular
\[
W(\theta(\cdot,z),\phi(\cdot,z))=1
\]
These functions are entire in $z$.  Note that $\phi(R,0)=\phi_0(R)$
from above. Furthermore,
\[
\theta_0(R):=\theta(R,0)= (1-2R^2+R^4/9)(1+R^2/3)^{-\frac32}
\]
The Weyl-Titchmarsh function $m(z)$ is uniquely defined by
\begin{equation}\label{eq:WT}
\psi_+(\cdot,z):=\theta(\cdot,z)+m(z) \phi(\cdot,z)\in L^2(0,\infty)
\;\forall \,\Im z>0
\end{equation}
The solution $\psi_+$ is referred to as the Weyl-Titchmarsh
solution. Then one has the following, see~\cite{GZ}:

\begin{prop}\label{prop:herglotz}
The function $m$ can be analytically continued to
$\Compl\setminus\spec(\cL)$ and it is a Herglotz function. For each
$R\ge0$, $\psi_+(R,z)$ and $\psi'_+(R,z)$ are analytic on
$\Compl\setminus\spec(\cL)$. The spectral measure of $\cL $ is
\[
d\rho = \delta_{\etao}+ \rho(\xi) d\xi, \qquad \rho(\xi)
=\frac1\pi\Im m(\xi+i0)
\]
in the following sense: the distorted Fourier transform defined as
  \begin{equation}\nonumber
    \calF: f\longrightarrow \hat{f}
 \end{equation}
  \begin{equation}\nonumber
    \hat{f}(\etao) =  \int_{0}^{\infty} \go(R) f(R)\,dR,
    \qquad      \hat{f} (\xi)=\lim_{b\rightarrow
      \infty}\int_{0}^{b}\phi(R,\xi)f(R)\,dR, \ \ \xi \geq 0
  \end{equation}
  is a unitary operator from $ L^2 (\R^+)$ to
  $L^{2}({\{\etao\}\cup {\R^+}},{\rho})$ and its inverse is given by
 \begin{equation}\nonumber
    \calF^{-1}: \hat{f}\longrightarrow f(R)=\hat{f}(\etao) \go(R)
    + \lim_{\mu \rightarrow
      \infty}\int_{0}^{\mu} \phi(R,
    \xi)\hat{f}(\xi)\,{\rho}(\xi)\,d\xi
  \end{equation}
Here $\lim$ refers to the $L^{2}({{\R^+}},{\rho})$, respectively the
$L^2(\R^+)$, limit.
\end{prop}

In the sequel we view the Fourier transform as a vector-valued map
\[
f\mapsto \binom{\hat{f}(\xi_d)}{\hat{f}(\cdot)}
\]
The Weyl-Titchmarsh solutions are scalar multiples of the Jost
solutions $f_+(R,z)$ which are determined via the condition that
\begin{equation}\label{eq:Jost}
  \cL f_+(\cdot,z)=z f_+(\cdot,z),\quad f_+(R,z)\sim e^{i\sqrt{z}R}
\text{\ \ as\ \ }R\to\infty
\end{equation}
where $\Im z\ge0, \Im\sqrt{z}\ge0$. They are solutions to the
integral equation, with $V(R)=-5(1+R^2/3)^{-2}$,
\[
f_+(R,z)=e^{i\sqrt{z}R} + \int_R^\infty
\frac{\sin(\sqrt{z}(R'-R))}{\sqrt{z}} V(R')f_+(R',z)\,dR'
\]
The functions $\psi_+(R,z)$ %and $f_+(R,z)$ have
 have well-defined limits as $z\to
\xi+i0$ when $\xi>0$. In particular, $\psi_+(R,\xi)\sim c_0(\xi)
e^{i\xi^{\frac12} R}$ as $R\to\infty$. The constant $c_0(\xi)$ is
determined from the Wronskian condition
$W(\psi_+(\cdot,\xi),\phi(\cdot,\xi))=1$. Once we have determined
$c_0(\xi)$ we find $m(\xi+i0)$ from the Wronskian relation
\[ m(\xi+i0) = W(\theta(\cdot,\xi), \psi_+(\cdot,\xi+i0)) \]

We now give an asymptotic expansion of our fundamental system for
small~$z$.

\begin{proposition}  \label{pphitheta} For any $z\in\Compl$
the fundamental system $\phi(R,z)$, $\theta(R,z)$ admits absolutely
convergent asymptotic expansions
\begin{align*}
\phi(R,z) &= \phi_0(R) + R^{-1} \sum_{j =1}^\infty (R^2
z)^{j} \phi_j(R^2)\\
 \theta(R,z) &= \theta_0(R) +  \sum_{j =1}^\infty (R^2
z)^{j} \theta_j(R^2)
\end{align*}
where the functions $\phi_j$, $\theta_j$ are holomorphic in $\Omega
= \{u\in\Compl: \Re u > -\frac12\} $ and satisfy the bounds
\begin{align}
| \phi_j(u)| &\leq \frac{C^j}{(j-1)!}|u|\la u\ra^{-\frac12},\label{eq:phij_bd}\\
% \quad |\phi_1(u)| > \frac{1}{2}u^{\frac12} \text{\ \ if\ \ }u\gg1 \\
 | \theta_j(u)| &\leq \frac{C^j}{(j-1)!} \la u\ra^{\frac12}, \qquad
u \in \Omega\label{eq:thetaj_bd}
\end{align}
Furthermore, \begin{equation}\label{eq:phi1} \phi_1(u) = \left\{
\begin{array}{ll} -\frac16 u (1+o(1))
&\text{\ \ as\ \ }u\to 0 \\
\frac{\sqrt{3}}{2} u^{\frac12}(1+o(1))
 &\text{\ \ as\ \ }u\to\infty \end{array} \right.
\end{equation}
\end{proposition}

\begin{proof}
We begin with $\phi$. We formally write
\[
\phi(R,z) = R^{-1} \sum_{j=0}^\infty z^j  f_j(R),\qquad
f_0(R)=R\phi_0(R)
\]
This becomes rigorous once we verify the convergence of the series
in any reasonable sense. The functions $f_j$ should solve
\[
\cL (R^{-1} f_j) = R^{-1} f_{j-1}, \quad  f_j(0)=f_j'(0)=0
\]
where we have set $f_{-1} = 0$.  The forward fundamental solution
for $\cL$ is
\[
H(R,R') = (\phi_0(R) \theta_0(R') -\phi_0(R') \theta_0(R))1_{[R
> R']}
\]
Hence we have the iterative relation
\[
f_j(R) = \int_0^R  \frac{R}{R'} \Big[\phi_0(R) \theta_0(R')
-\phi_0(R') \theta_0(R)\Big] f_{j-1}(R')\, dR', \quad f_0(R) =
R\phi_0(R)
\]
Using the expressions for $\phi_0$, $\theta_0$ we rewrite this as
\[\begin{split}
&f_j(R) =\\
& \int_0^R \frac{R^2(1-R^2/3)(1-2R'^2 +
R'^4/9)-RR'(1-R'^2/3)(1-2R^2+R^4/9)}{R'(1+R^2/3)^{\frac32}(1+R'^2/3)^{\frac32}}
f_{j-1}(R')\, dR'
\end{split}
\]
It is clear from this integral that each $f_j(R)$ is an analytic
function of $R^2$ provided $\Re R^2> -1$. Moreover, $f_j(R^2)$
vanishes like $R^{2(j+1)}$ around $R=0$, and grows like $R^{2j+1}$
as $R\to\infty$. The bounds in~\eqref{eq:phij_bd} are a quantitative
version thereof, and proved by induction.

For $\theta(R,z)$ we make the ansatz
\[
\theta(R,z) = \sum_{j=0}^\infty z^j  g_j(R),\quad g_0(R)=\theta_0(R)
\]
where the functions $g_j$ solve
\[
\cL g_j = g_{j-1}, \qquad g_{-1} = 0
\]
Thus,  the iterative relation is
\[
g_j(R) = \int_0^R  \Big[\phi_0(R) \theta_0(R') -\phi_0(R')
\theta_0(R)\Big] g_{j-1}(R')\, dR', \quad g_0(R) = \theta_0(R)
\]
The analyticity is the same as for the $f_j(R)$, and $g_j(R)$
vanishes like $R^{2j}$ around $R=0$, whereas the growth is
$R^{2j+1}$ as $R\to\infty$, as claimed.

Finally, the leading order of $\phi_1(u)$ is found by solving for
the coefficients $c_1 , c_2$ in
\begin{align*}
  (-\partial_{RR} - 5(1+R^2/3)^{-2}) c_1 R^3(1+o(1)) &= R(1+o(1))\text{\ \ as\ \ }R\to0 \\
(-\partial_{RR} - 5(1+R^2/3)^{-2}) c_2 R^2(1+o(1)) &=
-\sqrt{3}(1+o(1)) \text{\ \ as\ \ }R\to \infty,
\end{align*}
respectively. Thus, $c_1=-\frac16$, $c_2=\frac12\sqrt{3}$ as
claimed.
\end{proof}

Next, we turn to the asymptotic expansion of the Jost solutions.

\begin{proposition} For any $\xi>0$, the Jost solution  $f_+(\cdot,\xi)$
as in~\eqref{eq:Jost}  is of the form
\[
f_+(R,\xi) =  e^{iR \xi^\frac12} \sigma(R\xi^\frac12,R),\qquad
R^2\xi\gtrsim 1
\]
where $\sigma$ admits the asymptotic series approximation
\[
\sigma(q,R) \approx \sum_{j=0}^\infty q^{-j} \psi^+_j(R)
\]
in the sense that for all  integers $j_0\ge0$, and all indices
$\alpha$, $\beta$, we have
\begin{equation}\label{eq:symbol}
 \sup_{R>0}\,\la R\ra^{2}\Bigl|( R \partial_R)^\alpha (q
\partial_q)^\beta \Big[\sigma(q,R)
  - \sum_{j=0}^{j_0} q^{-j} \psi_j^+(R)\Big]\Bigr| \leq
c_{\alpha,\beta,j_0}\;   q^{-j_0-1}
\end{equation}
for all $q>1$. Here
\[
\psi^+_0 = 1, \quad \psi_1^+(R) = \left\{
\begin{array}{ll} ic_1R^{-2} + iO(R^{-4})
&\text{\ \ as\ \ }R\to\infty\\ ic_2R + iO(R^2)
 &\text{\ \ as\ \ }R\to 0 \end{array} \right.
\]
with some real constants $c_1,c_2$. More generally, $\psi^+_j(R)$
are smooth symbols of order $-2$ for $j\ge1$, i.e., for all $k\ge0$
\[\sup_{R>0} \la R\ra^2\,
|(\la R\ra \partial_R)^k \psi^+_j(R)| <\infty
\]
Finally, $\psi^+_j(R) = O(R^j)$ as $R\to0$.
 \label{ppsipsi}\end{proposition}

\begin{proof}
With the notation
\[
\sigma(q,R) =  f_+(R,\xi)  e^{-iR \xi^\frac12}
\]
we need to solve the conjugated equation
\begin{equation}
\left(-\partial_{RR} - 2 i \xi^\frac12 \partial_R  -
  \frac{5}{(1+R^2/3)^2}\right) \sigma(R\xi^{\frac12}, R) = 0
\label{conjug}\end{equation} We look for a formal power series
solving this equation,
\begin{equation}
\sum_{j=0}^\infty \xi^{-\frac{j}2} f_j(R) \label{formal}
\end{equation} This yields a recurrence relation for the $f_j$'s,
\[
2i \partial_R f_j = \left(-\partial_{RR}  -
  \frac{5}{(1+R^2/3)^2}\right) f_{j-1}, \quad
  f_j(\infty)=f_j'(\infty)=0
\]
with $f_0 = 1$,  which is solved by
\[
f_j(R) = \frac{i}{2} \partial_R f_{j-1}(R) -  \frac{i}{2}
\int_{R}^\infty
  \frac{5}{(1+R'^2/3)^2} f_{j-1}(R')\, dR'
\]
Then $f_1(R)$ is smooth for all $R\ge0$ with
$f_1(R)=-\frac{15i}{2R^3} + iO(R^{-5})$ as $R\to\infty$ and
$f_1(R)=ic+iO(R)$ around $R=0$. More generally, $f_j(R)= i^j
O(R^{-j-2})$ as $R\to\infty$ and $f_j(R)=i^j O(1)$ as $R\to0$.
Differentiating these functions leads to symbol-type behavior as
$R\to\infty$.  Defining $\psi_j^+(R):=R^j f_j(R)$ yields the desired
bounds.

To finish the proof, we first construct an approximate sum, i.e., a
function $\sigma_{ap}(q,R)$ with the property that for each $j_0
\geq 0$ we have
\begin{equation}\label{eq:sigma_schr}
\Big|(R \partial_R)^\alpha (q \partial_q)^\beta
\big[\sigma_{ap}(q,R)
  - \sum_{j=0}^{j_0} q^{-j} \psi_j^+(R)\big]\Big| \leq
c_{\alpha,\beta,j_0}\, \la R\ra^{-2} q^{-j_0-1}
\end{equation}
The construction of $\sigma_{ap}(q,R)$ is standard in symbol
calculus; we set
\[
\sigma_{ap}(q,R) := \sum_{j=0}^\infty q^{-j} \psi_j^+(R)
\chi(q\delta_j)
\]
where $\delta_j\to0$ sufficiently fast and $\chi$ is a cut-off
function which vanishes around zero and is equal to one for large
arguments. The bound~\eqref{eq:sigma_schr} implies that
$\sigma_{ap}(R\xi^\frac12,R)$ is a good approximate solution for
\eqref{conjug} at infinity, namely the error
\[
e(R \xi^\frac12,R) = \left(-\partial_{RR} -2 i \xi^\frac12
\partial_R
 - \frac{5}{(1+R^2/3)^2}\right) \sigma_{ap}(R,\xi)
\]
satisfies for all indices $\alpha,\beta,j$
\[
| (R \partial_R)^\alpha (q \partial_q)^\beta e(q,R) | \leq
c_{\alpha,\beta,j}\,  \la R\ra^{-4} q^{-j}
\]
To conclude the proof it remains to solve the equation for the
difference $\sigma_1 = -\sigma + \sigma_{ap}$,
\[
 \left(-\partial_{RR}- 2i \xi^\frac12 \partial_R
   - \frac{5}{(1+R^2/3)^2}\right) \sigma_1(R
\xi^\frac12,R) = e(R \xi^\frac12,R)
\]
with zero Cauchy data at infinity. We claim that the solution
$\sigma_1$ satisfies
\[
| ( R \partial_R)^\alpha (q \partial_q)^\beta \sigma_1(q,R) | \leq
c_{\alpha,\beta,j}\,  q^{-j}\la R\ra^{-2}, \qquad j \geq 2
\]
Note that this finishes the proof by defining $\sigma=\sigma_{ap}
-\sigma_1$.  A change of variable allows us to switch from the pair
of operators $(R\partial_R, q\partial_q)$ to $(R\partial_R, \xi
\partial_\xi)$ with comparable bounds. We rewrite the above equation as
a first order system for $(v_1,v_2) = (\sigma_1,R \partial_R
\sigma_1)$:
\[
\partial_R \left( \begin{array}{c} v_1 \cr v_2 \end{array}\right) -
\left( \begin{array}{cc} 0 & R^{-1} \cr
 - \frac{5R}{(1+R^2/3)^2}
 &R^{-1}  -2 i\xi^\frac12
 \end{array}\right)
 \left( \begin{array}{c} v_1 \cr v_2 \end{array}\right) =
 \left( \begin{array}{c} 0 \cr Re \end{array}\right)
\]
Then we have
\[
\frac{d}{dR} |v|^2 \gtrsim  -R^{-1} |v|^2 - R |v||e|
\]
which gives
\[
\frac{d}{dR} |v| \geq -C(R^{-1}|v| + R |e|)
\]
and by Gronwall
\[
|v(R)| \leq \int_{R}^\infty \left(\frac{R'}{R}\right)^C R' |e(R')|\,
dR'
\]
Then for large $j$ we have
\begin{equation}
|e| \lesssim \xi^{-\frac{j}{2}} R^{-j}\la R\ra^{-4} \implies |v|
\lesssim \xi^{-\frac{j}{2}} R^{-j}\la R\ra^{-2} \lesssim q^{-j} \la
R\ra^{-2}
 \label{ev}\end{equation}
 To estimate derivatives of $v$ we commute them with the operator. For
derivatives with respect to $R$ we have
\[\begin{split}
&\partial_R (R \partial_R) \left( \begin{array}{c} v_1 \cr v_2
\end{array}\right) - \left( \begin{array}{cc} 0 & \frac{1}R \cr
  - \frac{5R}{(1+R^2/3)^2}
 &\frac1R  -2 i\xi^\frac12
 \end{array}\right)
 (R \partial_R) \left( \begin{array}{c} v_1 \cr v_2
 \end{array}\right)\\
&= \left( \begin{array}{cc} 0 & \frac1R \cr
  \frac{10R(1-R^2/3)}{(1+R^2/3)^3}
 &\frac1R
 \end{array}\right)
 \left( \begin{array}{c} v_1 \cr v_2 \end{array}\right) +
 \left( \begin{array}{c} 0 \cr R \partial_R (Re) \end{array}\right)
\end{split}\] But the right-hand side is bounded by $R^{-j-1}$ from the
previous step and the hypothesis on $e$, therefore as above $R
\partial_R v$ is bounded by $R^{-j}$.

\noindent We argue similarly for the $\xi$ derivatives. We have
\[\begin{split}
&\partial_R (\xi \partial_\xi) \left( \begin{array}{c} v_1 \cr v_2
\end{array}\right) - \left( \begin{array}{cc} 0 & \frac{1}R \cr
 - \frac{5R}{(1+R^2/3)^2}
 &\frac1R  -2 i\xi^\frac12
 \end{array}\right)
 (\xi \partial_\xi) \left( \begin{array}{c} v_1 \cr v_2
 \end{array}\right)\\&
= \left( \begin{array}{cc} 0 & 0 \cr
 0 &   i\xi^\frac12
 \end{array}\right)
 \left( \begin{array}{c} v_1 \cr v_2 \end{array}\right) +
 \left( \begin{array}{c} 0 \cr \xi \partial_\xi (Re) \end{array}\right)
\end{split}\] The only difference is in the first term on the right, for
which we write $\xi^\frac12 = R^{-1} q$ and we use the decay
property of $v$ with $j$ replaced by $j+1$:
\[\begin{split}
|\xi^{\frac12} v_2| \less \xi^{\frac12} q^{-j-1}\less R^{-1}
q^{-j},\qquad  |\xi \partial_\xi (Re)| \less R^{-1} q^{-j}
\end{split}
\]
as desired.  Finally, higher order derivatives are estimated by
induction using the above arguments at each step.
\end{proof}

Next, we describe the spectral measure of $\cL$. Due to the
resonance at zero energy, the spectral density becomes unbounded for
small~$\xi$. In what follows, $f_-(R,\xi):= f_+(R,-\xi)=
\overline{f_+(R,\xi)}$.

\begin{lemma}\label{lem:spec_meas}
For all $\xi>0$ there is $a(\xi)\ne 0$ so that
\[
\phi(R,\xi) = a(\xi) f_+(R,\xi) + \overline{a(\xi)} f_-(R,\xi)
\]
with
\[
 |a(\xi)| \asymp \left\{\begin{array}{ll} 1
&\text{\ \ as\ \ }\xi\to 0 \\
\xi^{-\frac12} &\text{\ \ as\ \ }\xi \to\infty \end{array} \right.
\]
with symbol type behavior of all its derivatives. The density
$\rho(\xi)$ of the spectral measure satisfies
\[
\rho(\xi) \asymp \left\{\begin{array}{ll} \xi^{-\frac12}
&\text{\ \ as\ \ }\xi\to 0 \\
\xi^{\frac12} &\text{\ \ as\ \ }\xi \to\infty \end{array} \right.
\]
with symbol type behavior of all derivatives.
\end{lemma}
\begin{proof}
By inspection,
\[
a(\xi) = \frac{W(\phi(\cdot,\xi),
f_-(\cdot,\xi))}{W(f_+(\cdot,\xi),f_-(\cdot,\xi))} =
\frac{1}{-2i\xi^{\frac12}} W(\phi(\cdot,\xi), f_-(\cdot,\xi))
\]
By the preceding asymptotic analysis we can evaluate
\[
W(\phi(\cdot,\xi), f_+(\cdot,\xi)) =
\phi(\eps\xi^{-\frac12},\xi)f_+'(\eps\xi^{-\frac12},\xi) -
\phi'(\eps\xi^{-\frac12},\xi) f_+(\eps\xi^{-\frac12},\xi)
\]
with some small fixed $\eps>0$ to conclude that
\[
|W( \phi(\cdot,\xi), f_+(\cdot,\xi))| \less \left\{\begin{array}{ll}
\xi^{\frac12}
&\text{\ \ as\ \ }\xi\to 0 \\
1 &\text{\ \ as\ \ }\xi \to\infty \end{array} \right.
\]
with the corresponding upper bound on the derivatives. This yields
the desired upper bound on $|a(\xi)|$.  To obtain the lower bound,
we proceed as follows. First, observe that
\[
\Im (f_+(R,\xi) f_-'(R,\xi)) =-\xi^{\frac12}
\]
Second, it follows that
\[
\Im (f_-'(R,\xi)W (f_+(R,\xi),\phi(R,\xi))) =
\xi^{\frac12}\phi'(R,\xi)
\]
so that
\[
|W(f_+(R,\xi),\phi(R,\xi))| \ge \xi^{\frac12}
\frac{|\phi'(R,\xi)|}{|f_+'(R,\xi)|}
\]
>From our asymptotic analysis, again at $R=\eps \xi^{-\frac12}$,
\[
\frac{|\phi'(R,\xi)|}{|f_+'(R,\xi)|} \gtrsim
\left\{\begin{array}{ll} 1
&\text{\ \ as\ \ }\xi\to 0 \\
\xi^{-\frac12} &\text{\ \ as\ \ }\xi \to\infty \end{array} \right.
\]
which leads to the claimed lower bound on $|W( \phi(\cdot,\xi),
f_+(\cdot,\xi))| $.

\noindent The Weyl solution
\[
\psi_+(R,\xi+i0) = \theta(R,\xi) + m(\xi+i0)\phi(R,\xi)
\]
satisfies $\psi_+(\cdot,\xi+i0)=c_0(\xi) f_+(\cdot,\xi)$. Solving
for $c_0(\xi)$ yields
\[
m(\xi+i0)= \frac{W(\theta(\cdot,\xi),
f_+(\cdot,\xi))}{W(f_+(\cdot,\xi),\phi(\cdot,\xi))} =
\frac{W(\theta(\cdot,\xi), f_+(\cdot,\xi))W(
f_-(\cdot,\xi),\phi(\cdot,\xi))}{|W(f_+(\cdot,\xi),\phi(\cdot,\xi))|^2}
\]
Since
\[
f_+(\cdot,\xi) = -\phi(\cdot,\xi)
W(f_+(\cdot,\xi),\theta(\cdot,\xi)) + \theta(\cdot,\xi)
W(f_+(\cdot,\xi), \phi(\cdot,\xi))
\]
implies that
\[
-2i\xi^{\frac12}= W (f_+(\cdot,\xi), f_-(\cdot,\xi)) = -2i \Im [
W(\theta(\cdot,\xi), f_+(\cdot,\xi))W(
f_-(\cdot,\xi),\phi(\cdot,\xi))]
\]
we conclude that
\begin{equation}\label{eq:rho_m}
\rho(\xi) = \frac1\pi \Im m(\xi+i0) = \frac{\xi^{\frac12}}{\pi
|W(f_+(\cdot,\xi),\phi(\cdot,\xi))|^2} =
\frac{1}{4\pi}\frac{1}{\xi^{\frac12} |a(\xi)|^2}
\end{equation}
The denominator was estimated above, leading to the desired bound on
the spectral density.
\end{proof}

 \section{The transference identity}
 \label{sec:transference}

 Returning to the radiation part $\tilde{\epsilon}$ in
 \eqref{eq:lin_eq}, the idea is to expand it in terms of the
 generalized Fourier basis $\phi(R, \xi)$ from Proposition~\ref{prop:herglotz}, i.e., write
 \begin{equation}\nonumber
    \tileps(\tau, R)=x_0(\tau) \go(R) + \int_{0}^{\infty}  x(\tau,\xi)
   \phi(R,\xi)\rho(\xi)\,d\xi
 \end{equation}
 and deduce a transport equation for the Fourier coefficients
 $(x_0(\tau), x(\tau, \xi))$.  The main difficulty in doing this is
 caused by the operator $R \partial_R$ which is not diagonal in the
 Fourier basis. We re-express this derivative in terms of the
 derivative $2\xi\pr_\xi$. We refer to this procedure, which involves
 a certain error operator, as the transference identity since it
 allows us to transfer derivatives from $R$ to~$\xi$.  We define the
 error operator $\cK$ by
 \begin{equation}\label{eq:transfer}
 \widehat{R \partial_R u} = - 2 \xi \partial_\xi \hat u + \calK \hat u
 \end{equation}
 where $\hat{f}=\calF f$ is the ``distorted Fourier transform" from
 Proposition~\ref{prop:herglotz} and the operator $-2\xi \partial_\xi$
acts only on the continuous part of the spectrum.
 Apriori we have
 \[
 \calK: C_0^\infty(\{\etao\} \cup (0,\infty)) \to
 C^\infty( \{\etao\} \cup (0,\infty))
 \]
Splitting functions on $\spec(\cL)$ into a discrete and continuous
component we obtain a matrix representation for $\calK$,
\[
\calK = \left(\begin{array}{cc} \calK_{dd} & \calK_{dc} \cr
\calK_{cd} & \calK_{cc} \end{array} \right)
\]
Using the expressions for the direct and inverse Fourier transform
in Proposition~\ref{prop:herglotz} we obtain
\begin{align}
\calK_{dd}& =  \Bla R\partial_{R} \go(R), \go(R)\Bra_{L^2_R}\nn
\\
 \calK_{dc} f & =  \Bla \int_{0}^{\infty}f(\xi) R\partial_{R} \tilphi(R,
 \xi)\rho(\xi)\,d\xi\,,\, \go(R)\Bra_{L^2_R}\nn
\\
 \calK_{cd} (\eta)& = \Bla R\partial_{R} \go,  \tilphi(R,
 \eta) \Bra_{L^2_R}\nn
\\
 \calK_{cc} f(\eta) &=  \Bla
 \int_{0}^{\infty}f(\xi) R\partial_{R} \tilphi(R,
 \xi)\rho(\xi)\,d\xi\,,\, \tilphi(R, \eta)\Bra_{L^2_R}\nn\\
 &  + \Bla
 \int_{0}^{\infty}2\xi\partial_{\xi} f(\xi) \tilphi(R,
 \xi)\rho(\xi)\,d\xi\,,\, \tilphi(R, \eta)\Bra_{L^2_R}\label{eq:Kcc}
\end{align}
Integrating by parts with respect to $R$ in the first two relations we
obtain
\[
\calK_{dd} = -\frac12, \qquad  \calK_{dc} f = - \int_0^\infty f(\xi)
K_d(\xi) \rho(\xi) d\xi, \qquad  \calK_{cd} (\eta) = K_d(\eta)
\]
where
\[
K_d(\eta) =  \Bla R\partial_{R} \go,  \tilphi(R,
 \eta) \Bra_{L^2_R}
\]
Integrating by parts with respect to $\xi$ in \eqref{eq:Kcc} yields
 \begin{equation}\label{calk}\begin{split}
     \calK_{cc} f(\eta) &=  \Bla
     \int_{0}^{\infty}f(\xi)[R\partial_{R}-2\xi\partial_{\xi}]
\tilphi(R,\xi)\rho(\xi)\,d\xi\,,\,
     \tilphi(R, \eta)\Bra_{L^2_R}
\\ & - 2 \left(1+ \frac{\eta
         \rho'(\eta)}{\rho(\eta)}\right)  f(\eta)
   \end{split}
\end{equation}
 where the scalar product is to be interpreted in the principal value
 sense with $f\in C_0^\infty((0,\infty))$.

 In this section, we study the boundedness properties of the operator
 $\calK$. We begin with a description of the function $K_d$
and of the kernel $K_0(\eta,\xi)$ of $\calK_{cc}$.

 \begin{theorem}\label{tp}
a)   The operator $\calK_{cc}$ can be written as
   \begin{equation}\label{eq:kern}
     \calK_{cc} = -\Big(\frac{3}{2} + \frac{\eta \rho'(\eta)}{\rho(\eta)}\Big)\delta(\xi-\eta) + \calK_0
   \end{equation}
   where the operator $\calK_0$ has a kernel $K_0(\eta, \xi)$ of the
   form\footnote{The kernel below is interpreted in the principal
     value sense}
   \begin{equation} \label{ketaxi} K_0(\eta,
     \xi)=\frac{\rho(\xi)}{\eta-\xi} F(\xi,\eta)
   \end{equation}
   with a symmetric function $F(\xi,\eta)$ of class $C^2$ in
   $(0,\infty) \times (0,\infty)$ satisfying the bounds
   \[\begin{split}
   | F(\xi,\eta)| &\lesssim \left\{ \begin{array}{cc} \xi+\eta &
       \xi+\eta \leq 1 \cr (\xi+\eta)^{-1} (1+|\xi^\frac12
       -\eta^\frac12|)^{-N} & \xi+\eta \geq 1
     \end{array} \right.\\
   | \partial_{\xi} F(\xi,\eta)|+| \partial_{\eta} F(\xi,\eta)| &\lesssim \left\{
     \begin{array}{cc} 1 & \xi+\eta \leq 1 \cr (\xi+\eta)^{-\frac32}
       (1+|\xi^\frac12 -\eta^\frac12|)^{-N} & \xi+\eta \geq 1
     \end{array} \right.\\
  \sup_{j+k=2} | \partial^j_{\xi}\partial^k_{\eta} F(\xi,\eta)| &\lesssim \left\{
     \begin{array}{cc} (\xi+\eta)^{-\frac12} & \xi+\eta \leq 1 \cr
       (\xi+\eta)^{-2} (1+|\xi^\frac12 -\eta^\frac12|)^{-N} &
       \xi+\eta \geq 1
     \end{array} \right.
     \end{split}
   \]
   where $N$ is an arbitrary large integer.

b) The function $K_d$ is smooth and rapidly decaying at infinity.
\end{theorem}

 \begin{proof}
   We first establish the off-diagonal behavior of $\calK_{cc}$, and later
   return to the issue of identifying the $\delta$-measure that sits
   on the diagonal. We begin with \eqref{calk} with $f \in
   C_0^\infty((0,\infty))$. The integral
   \[
   u(R) =
   \int_{0}^{\infty}f(\xi)[R\partial_{R}-2\xi\partial_{\xi}]\tilphi(R,
   \xi)\rho(\xi)\,d\xi
   \]
   behaves like $R$ at $0$ and is a Schwartz function at
   infinity. The second factor $\tilphi(R,\eta)$ in \eqref{calk} also
   decays like $R$ at $0$ but at infinity it is only bounded
   with bounded derivatives. Then the following integration by parts
   is justified:
   \[
   \eta \calK_{cc} f(\eta) = \Bla u, \cL \tilphi(R, \eta)\Bra_{L^2_R} = \Bla
   \cL u, \tilphi(R, \eta)\Bra_{L^2_R}
   \]
   Moreover,
   \[
   \begin{split}
     \cL u =& \int_{0}^{\infty}f(\xi)[\cL,R\partial_{R}] \tilphi(R,
     \xi)\rho(\xi)\,d\xi +
     \int_{0}^{\infty}f(\xi)(R\partial_{R}-2\xi\partial_{\xi}) \xi
     \tilphi(R,\xi)\rho(\xi)\,d\xi \\ = & \int_{0}^{\infty}f(\xi) [\cL,R\partial_{R}]
     \tilphi(R, \xi)\rho(\xi)\,d\xi + \int_{0}^{\infty}\xi
     f(\xi)(R\partial_{R}-2\xi\partial_{\xi})
     \tilphi(R,\xi)\rho(\xi)\,d\xi \\
     & - 2\int_{0}^{\infty}\xi
     f(\xi)\tilphi(R,\xi)\rho(\xi)\,d\xi
   \end{split}
   \]
   with the commutator
   \[
   [\cL,R\partial_{R}] = 2\cL + \frac{10}{(1+R^2/3)^2} -
   \frac{20R^2}{3(1+R^2/3)^3} =: 2\cL + U(R)
   \]
    Thus,
    \[
\cL u =\int_{0}^{\infty}f(\xi) U(R)
     \tilphi(R, \xi)\rho(\xi)\,d\xi + \int_{0}^{\infty}\xi
     f(\xi)(R\partial_{R}-2\xi\partial_{\xi})
     \tilphi(R,\xi)\rho(\xi)\,d\xi
    \]
   Hence we obtain
   \[
   \eta \calK_{cc} f(\eta) - \calK_{cc} (\xi f)(\eta) = \Bla
   \int_{0}^{\infty}f(\xi) U(R) \tilphi(R, \xi)\rho(\xi)\,d\xi ,
   \tilphi(R, \eta)\Bra_{L^2_R}
   \]
   The double integral on the right-hand side is absolutely
   convergent, therefore we can change the order of integration to
   obtain
   \[
   (\eta -\xi) K_0(\eta,\xi) = \rho(\xi) \Bla U(R) \tilphi(R, \xi),
   \tilphi(R, \eta)\Bra_{L^2_R}
   \]
   This leads to the representation in \eqref{ketaxi} when
$\xi\ne\eta$ with
   \[
   F(\xi,\eta) = \Bla U(R) \tilphi(R, \xi), \tilphi(R,
   \eta)\Bra_{L^2_R}
   \]
   It remains to study its size and regularity. First, due
   to our pointwise bound from the previous section,
   \begin{equation}\begin{split}
   \sup_{R\ge0}|\tilphi(R, \xi)|
   &\less \la\xi\ra^{-\frac12},\\
    |R\partial_R \phi(R,\xi)|& \less R \qquad \forall\; \xi>1 \\
   |\partial_\xi \tilphi(R, \xi)|
   &\less \min(R\xi^{-1}, R^{3}) \qquad \forall\;
   \xi>1
   \\
|\partial_\xi \tilphi(R, \xi)|
   &\less \min(R\xi^{-\frac12},R^2) \qquad \forall\;
   0<\xi<1\\
|\partial^2_\xi \tilphi(R, \xi)|
   &\less \min(R^2\xi^{-\frac32}, R^{5}) \qquad \forall\;
   \xi>1
   \\
|\partial^2_\xi \tilphi(R, \xi)|
   &\less \min(R^2\xi^{-1}, R^4) \qquad \forall\;
   0<\xi<1
   \end{split}\label{eq:phi_est}
   \end{equation}
   we always have the estimates
   \begin{equation}\label{eq:Fbd}\begin{split}
   |F(\xi,\eta)| &\less \la\xi\ra^{-\frac12}\la\eta\ra^{-\frac12},\\
 |\partial_\xi F(\xi,\eta)| &\less
 \la\xi\ra^{-1}\la\eta\ra^{-\frac12}, \quad
 |\partial_\eta F(\xi,\eta)| \less
 \la\xi\ra^{-\frac12}\la\eta\ra^{-1}, \\
 |\partial_{\xi\eta} F(\xi,\eta)| &\less
\xi^{-1}\eta^{-1} \qquad \forall\;
 \xi>1,\,\eta>1\\
|\partial_\xi^2 F(\xi,\eta)| &\less \xi^{-\frac32}\eta^{-\frac12} \qquad \forall\;\xi>1,\,\eta>1\\
|\partial_\eta^2 F(\xi,\eta)| &\less \xi^{-\frac12}\eta^{-\frac32}
\qquad \forall\;\xi>1,\,\eta>1
\end{split}
\end{equation}
   They are only useful when $\xi$ and $\eta$ are very close. To improve on them, we consider two
   cases:

   {\bf Case 1: $1 \lesssim \xi+\eta $.}
   To capture the cancelations when $\xi$ and $\eta$ are separated we
   resort to another integration by parts,
   \begin{equation}
   \label{eq:Fdef}
   \eta F(\xi,\eta) = \Bla U(R) \tilphi(R, \xi), \cL \tilphi(R,
   \eta)\Bra = \Bla [\cL, U(R)] \tilphi(R, \xi), \tilphi(R, \eta)\Bra +
   \xi F(\xi,\eta)
   \end{equation}
   Hence, evaluating the commutator,
   \begin{equation} \label{Kk1} (\eta -\xi) F(\xi,\eta) = -\Bla (2 U_R
     \partial_R + U_{RR}) \tilphi(R, \xi), \tilphi(R, \eta)\Bra
   \end{equation}
   Since $U_R(0)=0$ it follows that $ (2 U_R \partial_R + U_{RR})
   \tilphi(R, \xi)$ has the same behavior as $\tilphi(R, \xi)$ in the
   first region. Then we can repeat the argument above to obtain
   \[
   (\eta -\xi)^2 F(\xi,\eta) = -\Bla [\cL,2 U_R \partial_R + U_{RR}]
   \tilphi(R, \xi), \tilphi(R, \eta)\Bra
   \]
   The second commutator has the form, with $V(R):= -
   5(1+R^2/3)^{-2}$,
   \[ [\cL,2 U_R \partial_R + U_{RR}] = 4 U_{RR} \cL - 4U_{RRR}
   \partial_R -U_{RRRR}
    - 2 U_R V_R -4U_{RR} V
   \]
   Since $V(R),U(R)$ are even,  this leads to
   \[
   (\eta -\xi)^2 F(\xi,\eta) = \Bla ( U^{odd}(R) \partial_R + U^{even}(R) +
   \xi U^{even}(R) ) \tilphi(R, \xi), \tilphi(R, \eta)\Bra
   \]
   where by $ U^{odd}$, respectively $U^{even}$, we have generically
   denoted odd, respectively even, nonsingular rational functions with
   good decay at infinity.  Inductively, one now verifies the identity
        \begin{equation} \begin{split}
     &(\eta -\xi)^{2k} F(\xi,\eta) =  \Bla \Big( \sum_{j=0}^{k-1} \xi^{j}\, U_{kj}^{odd}\, \partial_R +
     \sum_{\ell=0}^k \xi^\ell U_{k\ell}^{even} \Big) \tilphi(\cdot, \xi), \tilphi(\cdot, \eta)\Bra
     \label{Kk} \\
&\la R\ra|U_{kj}^{odd}(R)| + |U_{k\ell}^{even}(R)| \less \la
R\ra^{-4-2k} \qquad \forall\; j,\ell
     \end{split}\end{equation}
   By means of the pointwise bounds on $\tilphi$ and $\partial_R
   \tilphi$ from~\eqref{eq:phi_est} we infer from this that
   \[
   |F(\xi,\eta)| \lesssim \frac{\la\xi\ra^{k -\frac12}
     \la\eta\ra^{-\frac12}} {(\eta -\xi)^{2k}} \qquad\forall \;
     \xi, \;\eta>0
   \]
Combining this estimate with \eqref{eq:Fbd} yields, for arbitrary
$N$,
\[ |F(\xi,\eta)|\less (\xi+\eta)^{-1} (1+|\xi^\frac12
       -\eta^\frac12|)^{-N} \text{\ \ provided\ \ } \xi+\eta \gtrsim 1,
       \]
as claimed.
   For the derivatives of $F$ we follow a similar procedure. If $\xi$
   and $\eta$ are comparable,  then from~\eqref{eq:Fbd},
   $|\partial_\eta F(\xi,\eta)|\less \la \xi\ra^{-\frac32}$.
   Otherwise we differentiate with respect to $\eta$ in \eqref{Kk}.
   This yields
   \[\begin{split}
   (\eta -\xi)^{2k} \partial_\eta F(\xi,\eta) &=
   \Bla \Big( \sum_{j=0}^{k-1} \xi^{j}\, U_{kj}^{odd}\, \partial_R +
     \sum_{\ell=0}^k \xi^\ell U_{k\ell}^{even} \Big)\tilphi(R, \xi),
   \partial_\eta \tilphi(R, \eta)\Bra \\
   &\quad - 2k (\eta -\xi)^{2k-1}
   F(\xi,\eta)
   \end{split}
   \]
   Using also the bound on $F$ from above we obtain
   \[
   |\partial_\eta F(\xi,\eta)| \lesssim \frac{\xi^{k-\frac12}\eta^{-1}+\xi^{k-1}\eta^{-\frac12}} {(\eta -\xi)^{2k}}, \qquad 1
   \lesssim \xi,\eta
   \]
   respectively
   \[
   |\partial_\eta F(\xi,\eta)| \lesssim\frac{
     \eta^{-1}} {(\eta -\xi)^{2k}} \qquad \xi \ll 1 \lesssim
   \eta
   \]
   and
   \[
   |\partial_\eta F(\xi,\eta)| \lesssim\frac{\xi^{k-\frac12}}{(\eta -\xi)^{2k}} \qquad \eta \ll 1 \lesssim
   \xi
   \]
   which again yield the desired bounds.
   Finally, we consider the second order derivatives with respect to
   $\xi$ and $\eta$. For $\xi$ and $\eta$ close we again use the
   bound from~\eqref{eq:Fbd}.
   Otherwise we differentiate twice in \eqref{Kk} and continue as
   before.  We note that it is important here that the decay of
   $U_{kj}^{odd}$ and $U_{k\ell}^{even}$ improves with $k$. This is because the
   second order derivative bound at $0$ has a sizeable growth at
   infinity which has to be canceled,
   \[
   | \partial_\xi^2 \tilphi(R, 0)| \approx R^4
   \]

   {\bf Case 2: $ \xi,\eta \ll 1$}. First, we note that $F(0,0) =
   0$. This can be verified by direct integration, and is
   heuristically justified by the fact that $U = [\cL,R\partial_R]-2\cL$.
   The pointwise bound
   \[
   \sup_{\xi,\eta>0}| \partial_\xi F(\xi,\eta)| \lesssim 1
   \]
   follows by differentiating~\eqref{eq:Fdef} and from the bound
   $|\pr_\xi\phi(R,\xi)|\less R^2$, see~\eqref{eq:phi_est}.
   To bound the second order derivatives of $F$ we recall the pointwise
   bounds, for $0<\xi<1$,
   \[
   |\partial^j_\xi \tilphi(R, \xi)| \lesssim \left\{ \begin{array}{cc}
       R^{2j}  & R < \xi^{-\frac12} \\ \xi^{-j/2}
       \, R^j & R \geq \xi^{-\frac12}
     \end{array} \right., \qquad j = 0,1,2
   \]
   If $\eta < \xi<1$, then these bounds imply that
   \[\begin{split}
   | \partial_{\xi\eta} F(\xi,\eta)| &\lesssim
   \int_{0}^{\xi^{-\frac12}} \la R\ra^{-4} R^4\, dR +
   \int_{\xi^{-\frac12}}^{\eta^{-\frac12}} \la R\ra^{-4} R^{3} \xi^{-\frac12}
   \,dR
   +\int_{\eta^{-\frac12}}^\infty \la R\ra^{-2} (\xi\eta)^{-\frac12}
   \, dR \\
   & \less [1+\log(\xi/\eta)] \xi^{-\frac12}
   \end{split}
   \]
   The logarithm in the middle integral is an artefact and can be
   removed using the oscillatory nature of $\partial_\xi\phi(R,\xi)$
   in the regime $R^2\xi>1$ as provided by Proposition~\ref{ppsipsi}
   and Lemma~\ref{lem:spec_meas}. Loosely speaking, this means
   integrating by parts using
   that $\partial_\xi \phi(R,\xi)\sim R\xi^{-1}\partial_R
   e^{iR\xi^{\frac12}}$ for $R^2\xi>1$ and small $\xi$.
   Thus, actually
   \[
| \partial_{\xi\eta} F(\xi,\eta)| \less \xi^{-\frac12}
   \]
   A similar computation yields, for $\xi<1$,
\[\begin{split}
   | \partial^2_{\xi} F(\xi,\eta)| &\lesssim
   \int_{0}^{\xi^{-\frac12}} \la R\ra^{-4} R^4 \, dR +
   \int_{\xi^{-\frac12}}^{\infty} \la R\ra^{-4} R^2\xi^{-1}
   \,dR\less \xi^{-\frac12}
   \end{split}
   \]
   This bound is too weak when $\xi\ll \eta<1$.
   In that case, we  differentiate \eqref{Kk1} to
   obtain
   \[
   (\eta-\xi) \partial_\xi^2 F(\xi,\eta) = 2 \partial_\xi
   F(\xi,\eta) + \Bla \partial_\xi^2 \tilphi(R, \xi), (2 U_R \partial_R + U_{RR})
    \tilphi(R, \eta)\Bra
   \]
  which in turn yields
   \begin{equation}\label{eq:letzt}
   (\eta-\xi) \partial_\xi^2 F(\xi,\eta) = \int_\xi^\eta \Big[2
   \partial_{\xi\zeta} F(\xi,\zeta) + \Bla  \partial^2_\xi \tilphi(\cdot, \xi),
   (2 U_R \partial_R
   + U_{RR}) \partial_\zeta
   \tilphi(\cdot, \zeta)\Bra\Big] \, d\zeta
   \end{equation}
   Using also the bound
   \[
   | R\partial_{R\zeta} \tilphi(R,\zeta) | \less
   \min(R\zeta^{-\frac12},R^2)
   \]
   we can evaluate the inner product in \eqref{eq:letzt} as follows:
\begin{align*}
&\Big|\Bla \partial_\xi^2 \tilphi(\cdot, \xi), (2 U_R \partial_R +
U_{RR}) \pr_\zeta\tilphi(\cdot, \zeta)\Bra\Big| \\
&\less \int_0^{\zeta^{-\frac12}} \la R\ra^{-6} R^4 R^2 \, dR +
\int_{\zeta^{-\frac12}}^{\xi^{-\frac12}} \la R\ra^{-6} R^4
R\zeta^{-\frac12} \, dR + \int_{\xi^{-\frac12}}^\infty \la R\ra^{-6}
R^2\xi^{-1} R\zeta^{-\frac12} \, dR\\
&\less [1+\log(\zeta/\xi)]\zeta^{-\frac12}
\end{align*}
The logarithm appearing in the  middle integral can be removed as
before exploiting cancellations. Thus, \eqref{eq:letzt} is
controlled by
   \[
   |(\eta-\xi) \partial_\xi^2 F(\xi,\eta)| \lesssim \Big|\int_\xi^\eta
    \zeta^{-\frac12} \,d\zeta\Big| \lesssim
    \eta^{\frac12}
   \]
   Since $\eta \gg \xi$ this yields
   \[
   |\partial_\xi^2 F(\xi,\eta)| \lesssim \eta^{-\frac12}
   \]
   which concludes the analysis of the off-diagonal part of the kernel.

   Next, we extract the $\delta$ measure that sits on the diagonal of the
   kernel of $\calK_{cc}$ from the representation formula~\eqref{calk}, see also~\eqref{eq:kern}.
   To do so, we can restrict
    $\xi,\eta$  to a compact subset of
   $(0,\infty)$.  This is convenient, as we then have the following
   asymptotics of $\tilphi(R, \xi)$ for $R \xi^\frac12 \gg 1$:
   \[ \begin{split}
   \tilphi(R,\xi) &=
   2\Re [a(\xi)
     e^{iR\xi^{\frac{1}{2}}} ] + O(R^{-2})\\
   (R\partial_R -2\xi \partial_\xi) \tilphi(R,\xi) &= -4\Re [ \xi
     a'(\xi)
     e^{i R\xi^{\frac{1}{2}}} ] + O(R^{-2})
       \end{split}
   \]
   where the $O(\cdot)$ terms depend on the choice of the compact
   subset.
   The $R^{-2}$ terms are integrable so they contribute a bounded
   kernel to the inner product in~\eqref{calk}. The same applies to the
contribution of a bounded $R$ region.  Using the above expansions,
we conclude that the $\delta$-measure contribution of the inner
   product in~\eqref{calk} can only come from one of the following
   integrals:
   \begin{align}
     &  -4\int_0^\infty \int_0^\infty f(\xi) \chi(R) \Re
     \big[ \xi {a}'(\xi)
       a(\eta)
       e^{iR(\xi^{\frac{1}{2}}+\eta^\frac12)}\big]\rho(\xi)\, d\xi dR \label{eq:I1}\\
       &  -4\int_0^\infty \int_0^\infty f(\xi) \chi(R)\,
      \Re\big[ \xi {a}'(\xi)
       \bar{a}(\eta)
       e^{iR(\xi^{\frac{1}{2}}-\eta^\frac12)}\big]\;\rho(\xi)\, d\xi dR \label{eq:I2}
   \end{align}
   where $\chi$ is a smooth cutoff function which equals $0$ near
   $R=0$ and $1$ near $R=\infty$.  In all of the above integrals we can
   argue as in the proof of the classical Fourier inversion formula to
   change the order of integration. Integrating by parts in the first
   integral~\eqref{eq:I1} reveals that it cannot
   contribute a $\delta$-measure. On the other hand, \eqref{eq:I2}
contributes both a Hilbert transform type kernel  as well as a
$\delta$-measure to $K$. By inspection, the $\delta$ contribution is
\[\begin{split}
&-2 \int_{-\infty}^\infty \Re
     [ \xi {a}'(\xi)
       \bar{a}(\eta)
       e^{iR(\xi^{\frac{1}{2}}-\eta^\frac12)}]\rho(\xi)\, dR\\
      &= -4\pi \Re
     [ \xi {a}'(\xi)
       \bar{a}(\eta)]\rho(\xi)\delta(\xi^{\frac12}-\eta^{\frac12})\\
       &= -8\pi \xi^{\frac12}  \Re
     [ \xi{a}'(\xi)
       \bar{a}(\xi)
       ]\rho(\xi)\delta(\xi-\eta)\\
       &= \Big[\frac12 +
       \frac{\xi\rho'(\xi)}{\rho(\xi)}\Big]\delta(\xi-\eta)
\end{split}\]
where we used that $\rho(\xi)^{-1} = 4\pi \xi^{\frac12}|a(\xi)|^{2}$
in the final step, see~\eqref{eq:rho_m}. Combining this with the
$\delta$-measure in \eqref{calk} yields~\eqref{eq:kern}.

b)  Arguing as in part (a) we have
\[
K_d(\eta) = \frac{F(\etao,\eta)}{\etao-\eta}
\]
For $F$ we use the representation in \eqref{Kk} with $\xi$ replaced
by $\etao$ and $\phi(\cdot,\xi)$ replaced by $\phi_d$.
The conclusion easily follows from pointwise bounds on
$\phi(\cdot,\eta)$ and its derivatives.
 \end{proof}

 Next we consider the $L^2$ mapping properties for $\calK$. We
 introduce the weighted $L^2$ spaces $L^{2,\alpha}_\rho$
of functions on $\spec(\cL)$  with norm
 \begin{equation}\label{eq:L2weight}
 \| f\|_{L^{2,\alpha}_\rho}^2 :=  |f(\etao)|^2 +\int_0^\infty |f(\xi)|^2
 \la \xi\ra^{2\alpha} \rho(\xi)\,d\xi
 \end{equation}
 Then we have

\begin{proposition}
 a)  The operators $\calK_0$, $\calK$  map
  \[
  \calK_0\::\: L^{2,\alpha}_\rho \to L^{2,\alpha+1/2}_\rho, \qquad
 \calK \::\: L^{2,\alpha}_\rho \to L^{2,\alpha}_\rho.
  \]
b) In addition, we have the commutator bound
\[
[ \calK,\xi \partial_\xi] \::\:  L^{2,\alpha}_\rho \to
L^{2,\alpha}_\rho
\]
with $\xi\partial_\xi$ acting only on the continuous spectrum. Both
statements hold for all $\alpha\in\R$. \label{l2k}\end{proposition}
\begin{proof} We commence with the $\calK_0$ part.
a) The first property is equivalent to showing that the kernel
\[
\rho^{\frac12}(\eta)\la \eta\ra^{\alpha+1/2} K_0(\eta,\xi) \la
\xi\ra^{-\alpha} \rho^{-\frac12}(\xi) \::\: L^2(\R^+)\to L^2(\R^+)
\]
With the notation of the previous theorem, the kernel on the
left-hand side is
\[
\tilde K_0(\eta,\xi):= \la\eta\ra^{\alpha+1/2}\la\xi\ra^{-\alpha}
\frac{\sqrt{\rho(\xi)\rho(\eta)}}{\xi-\eta} F(\xi,\eta)
\]
We first separate the diagonal and off-diagonal behavior of $\tilde
K_0$, considering several cases.

  { \bf Case 1: $(\xi,\eta) \in Q := [0,4] \times [0,4]$.}

  We cover the unit interval with dyadic subintervals $I_j =
  [2^{j-1},2^{j+1}]$. We cover the diagonal with the union of squares
  \[
  A = \bigcup_{j=-\infty}^2 I_j \times I_j
  \]
  and divide the kernel $\tilde K_0$ into
  \[
  1_Q\tilde K_0 = 1_{A\cap Q} \tilde K_0 + 1_{Q \setminus A} \tilde K_0
  \]

  { \bf Case 1(a)}: Here we show that the diagonal part $1_{A\cap Q} \tilde K_0$ of
  $\tilde K_0$ maps $L^2$ to $L^2$. By orthogonality it suffices to restrict
  ourselves to a single square $I_j \times I_j$. We recall the $T1$
  theorem for Calderon-Zygmund operators, see page~293 in~\cite{Stein}: suppose the kernel
  $K(\eta,\xi)$ on $\R^2$ defines an operator $T:\calS\to\calS'$ and has the
  following pointwise properties with some $\gamma\in(0,1]$ and a constant $C_0$:
  \begin{enumerate}
  \item[(i)] $|K(\eta,\xi)|\le C_0|\xi-\eta|^{-1}$
  \item[(ii)] $|K(\eta,\xi)-K(\eta',\xi)|\le
  C_0|\eta-\eta'|^\gamma|\xi-\eta|^{-1-\gamma}$ for all
  $|\eta-\eta'|<|\xi-\eta|/2$
\item[(iii)] $|K(\eta,\xi)-K(\eta,\xi')|\le
  C_0|\xi-\xi'|^\gamma|\xi-\eta|^{-1-\gamma}$ for all
  $|\xi-\xi'|<|\xi-\eta|/2$
\end{enumerate}
If in addition $T$ has the restricted $L^2$ boundedness property,
i.e., for all $r>0$ and $\xi_0,\eta_0\in\R$,
$\|T(\omega^{r,\xi_0})\|_2 \le C_0r^{\frac12}$ and
$\|T^*(\omega^{r,\eta_0})\|_2 \le C_0r^{\frac12}$ where
$\omega^{r,\xi_0}(\xi)= \omega((\xi-\xi_0)/r)$ with a fixed
bump-function $\omega$, then $T$ and $T^*$ are $L^2(\R)$ bounded
with an operator norm that only depends on~$C_0$.

Within the
  square $I_j\times I_j$,
  Theorem~\ref{tp} shows that the kernel of $\tilde K_0$ satisfies these properties with $\gamma=1$,
  and is thus bounded on $L^2$.

  { \bf Case 1(b)}: Consider now the off-diagonal part $1_{Q \setminus
    A} \tilde K_0$. In this region, by Theorem~\ref{tp},
    \[
|\tilde K_0(\eta,\xi)| \less (\xi\eta)^{-\frac14}
    \]
which is a Hilbert-Schmidt kernel on $Q$ and thus $L^2$ bounded.

  {\bf Case 2:} $(\xi,\eta) \in Q^c $. We cover the diagonal with the union of squares
  \[
  B = \bigcup_{j=1}^\infty I_j \times I_j
  \]
  and divide the kernel $\tilde K_0$ into
  \[
  1_{Q^c}\tilde K_0 = 1_{B\cap Q^c} \tilde K_0 + 1_{Q^c \setminus B} \tilde K_0
  \]

  {\bf Case 2a:} Here we consider the estimate on $B$. As in case 1a)
  above, we use Calderon-Zygmund theory. Evidently, $|\tilde
  K_0(\eta,\xi)|\less |\xi-\eta|^{-1}$ on $B$ by Theorem~\ref{tp}.
  To check (ii) and (iii), we differentiate $\tilde K_0$. It will suffice to
  consider the case where the $\partial_\xi$ derivative falls on $F(\xi,\eta)$.
  We distinguish two cases:
if $|\xi^{\frac12}-\eta^{\frac12}|\le 1$, then $|\xi-\eta|\less
\xi^{\frac12}$ which implies that
\[
\frac{\xi^{-\frac12}|\xi-\xi'|}{|\xi-\eta|} \less
\frac{|\xi-\xi'|^{\frac12}}{|\xi-\eta|^{\frac32}} \qquad \forall\;
|\xi-\xi'|<|\xi-\eta|/2
\]
if, on the other hand, $|\xi^{\frac12}-\eta^{\frac12}|> 1$, then
\[
\frac{\xi^{-\frac12}|\xi-\xi'|}{|\xi-\eta||\xi^{\frac12}-\eta^{\frac12}|}
\less \frac{|\xi-\xi'|}{|\xi-\eta|^2} \qquad \forall\;
|\xi-\xi'|<|\xi-\eta|/2
\]
which proves property (iii) on $B$ with $\gamma=\frac12$, and by
symmetry also (ii). The restricted $L^2$ property follows form the
cancelation in the kernel and the previous bounds on the kernel.
Hence, $\tilde K_0$ is $L^2$ bounded on~$B$.

  {\bf Case 2b:} Finally, in the exterior region $Q^c \setminus B$ we
  have the bound, with arbitrarily large $N$,
  \[
  |\tilde K_0(\eta,\xi)| \lesssim  (1+\xi)^{-N} (1+\eta)^{-N}
  \]
  which is $L^2$ bounded by Schur's lemma.

This concludes the proof of the first mapping property in
part (a). The second one follows in a straightforward manner
since $K_d$ is rapidly decaying at $\infty$.

b) A direct computation shows that the kernel $K_0^{com}$ of the
commutator $[\xi \partial_\xi, K_0]$  is given by
\[
K_0^{com}(\eta,\xi) = (\eta \partial_\eta + \xi \partial_\xi)
K_0(\eta,\xi) + K_0(\eta,\xi) =
\frac{\rho(\xi)}{\xi-\eta}F^{com}(\xi,\eta)
\]
interpreted in the principal value sense and with $F^{com}$ given by
\[
F^{com}(\xi,\eta) = \frac{\xi \rho'(\xi)}{\rho(\xi)} F(\xi,\eta) +
(\xi \partial_\xi+\eta
\partial_\eta) F(\xi,\eta)
\]
By Theorem~\ref{tp} this satisfies the same pointwise off-diagonal
bounds as $F$. Near the diagonal the bounds for $F^{com}$ and its
derivatives are worse\footnote{The one derivative loss can be avoided
  by a more careful analysis, but this does not seem necessary here.}
than those for $F$ by a factor of $(1+\xi)^\frac12$. Then the proof of
the $L^2$ commutator bound for $K_0$ is similar to the argument in
part~(a).

The remaining part of the commutator $[\calK,\xi \partial_\xi]$
involves

(i) The commutator of the diagonal part of $\calK_{cc}$ with $\xi
\partial_\xi$. This is the multiplication operator by
\[
\xi \partial_\xi \frac{\xi \rho'(\xi)}{\rho(\xi)}
\]
which is bounded since $\rho$ has symbol like behavior both at $0$
and at $\infty$.

(ii) The operator $\xi \partial_\xi \calK_{cd}$ which is given by the
bounded rapidly decreasing function $\xi \partial_\xi K_d(\xi)$.

(iii)  The operator $\calK_{dc}\xi \partial_\xi $ given by
\[
\calK_{dc}\xi \partial_\xi f = \int_0^\infty K_d(\xi) \xi \partial_\xi f(\xi) d\xi
= - \int_0^\infty  f(\xi) \partial_\xi (\xi K_d(\xi))d\xi
\]
which is also bounded due to the properties of $K_d$.
\end{proof}

\section{The final system of equations}

We now rewrite the main equation \eqref{eq:lin_eq} in terms of the
Fourier transform. With $\calF$ as in Proposition~\ref{prop:herglotz},
and with $\beta=\dot\lambda\lambda^{-1}$,
\[
 \calF \Big(\partial_{\tau}+\beta(\tau) R\partial_{R}\Big) =
\Big(\partial_{\tau} - 2 \beta(\tau)\xi
\partial_{\xi}
 \Big) \calF  + \beta \calK  \calF
\]
which gives
\[
\begin{split}
&\calF \Big(\partial_{\tau}+\beta R\partial_{R}\Big)^2 =
 \Big(\partial_{\tau} + \beta (-2\xi \partial_{\xi}
 + \calK) \Big)^2 \calF
\\ = &  \Big(\partial_{\tau} - 2\beta\xi
\partial_{\xi}\Big)^2  \calF +
2\beta \calK  \Big(\partial_{\tau} - 2\beta\xi
\partial_{\xi}\Big)  \calF +
\beta^2(\calK^2 + 2[\xi
\partial_\xi,\calK]) \calF
\end{split}
\]
Recall that
\begin{equation}\nonumber
    \tileps(\tau, R)=x_0(\tau) \go + \int_{0}^{\infty}  x(\tau,\xi)
   \phi(R,\xi)\rho(\xi)\,d\xi
 \end{equation}
 This leads to a transport type equation for the Fourier
transform
\[
X(\tau) = (x_0(\tau), x(\tau,\xi))
\]
 of $ \tileps$ by applying $\calF$
to~\eqref{eq:lin_eq}. It is convenient to write it as a system for
the two components:
\begin{equation}\label{final}
\begin{split}
&\left( \begin{array}{cc} \partial_\tau^2 +\etao& 0 \cr 0 &
\Big(\partial_{\tau} - 2\beta\xi
\partial_{\xi}\Big)^2 +\xi \end{array}
\right) X\\ &= \beta(I- 2\calK)  \Big(\partial_{\tau} - 2\beta\xi
\partial_{\xi}\Big) X - \beta^2(\calK^2 -\calK +
2[\xi \partial_\xi,\calK]) X \\
&
 \quad + \lambda^{-2} \calF R (N_{2k-1} (R^{-1}
\calF^{-1} X) +  e_{2k-1})
\end{split}
\end{equation}
where it is understood that
\[
N_{2k-1} (R^{-1} \calF^{-1} X) +  e_{2k-1} = (N_{2k-1} (R^{-1}
\calF^{-1} X) +  e_{2k-1})(t(\tau),\lambda^{-1}(\tau)R)
\]
Note that $N_{2k-1}$ and $e_{2k-1}$ are only defined on $R\less
\tau$, but for the Fourier transform we need to extend them to all
$R$ -- this will be described in the next section, but for the
moment just take an arbitrary compactly supported extension with
the same regularity.

 We treat
problem~\eqref{final} iteratively, as a small perturbation of the
linear equation governed by the operator on the left--hand side. For
this we need to solve the following uncoupled system consisting of
an {\bf elliptic equation } and a {\bf{transport equation}}:
\begin{equation}\label{transport}
\begin{split}
 \Big[ \partial_{\tau}^2+\etao \Big]x_d(\tau)=b_d(\tau), \\
 \Big[\Big(\partial_{\tau}-
2\beta(\tau)\xi\partial_{\xi}\Big)^{2}+\xi\Big]x(\tau, \xi)=b(\tau,
\xi),
\end{split}
\end{equation}
We want to obtain solutions to \eqref{final} which decay as $\tau
\to \infty$. For the first equation above this is achieved by using
the standard fundamental solution $H_0$ which has kernel
\[
H_0(\tau,s) = -\frac12 |\etao|^{-\frac12}
e^{-|\etao|^\frac12|\tau-\sigma|}
\]
This means that up to homogeneous solutions of the form
$e^{-|\xi_d|^{\frac12}\tau}$ the unique bounded solution to the
elliptic equation is
\[
x_d(\tau)=-\frac12 |\etao|^{-\frac12} \int_0^\infty
e^{-|\etao|^\frac12|\tau-\sigma|}\; b_d(\sigma)\, d\sigma
\]
 The second equation is restricted to the range
$\xi
> 0$. Hence it is hyperbolic in nature, which means we can solve it
backward in time, i.e., with zero Cauchy data at $\tau = \infty$. We
denote by $H$ the backward fundamental solution for the operator
\[
 \Big(\partial_{\tau} - 2\beta(\tau)\xi
\partial_{\xi}\Big)^2   +\xi
\]
and by $H(\tau,\sigma)$ its kernel,
\[
x(\tau) = \int_{\tau}^\infty H(\tau,\sigma) f(\sigma)\, d\sigma
\]
Combining the two components we obtain a fundamental solution for the
system,
\[
\calH = {\rm diag}(H_0,H)
\]
Then we look for a solution $X$ to \eqref{final} as a solution to
the fixed point problem
\begin{equation}\label{final1}
\begin{split}
X &= \calH\Big( \beta(I- 2\calK)  \Big(\partial_{\tau} - 2\beta \xi
\partial_{\xi}\Big) X - \beta^2(\calK^2 -\calK +
2[\xi \partial_\xi,\calK]) X \\
&  + \lambda^{-2} \calF R (N_{2k-1} (R^{-1} \calF^{-1} X) +
e_{2k-1}) \Big)
\end{split}
\end{equation}

\begin{remark}
One can also freely add $ C e^{-|\etao|^\frac12\tau}$ to the first
component $x_d$ of $X$. Thus the fixed point argument yields in effect
a one parameter family of solutions $X$ depending on the parameter $C$.
\end{remark}

The mapping properties of $H$ are described in the following result,
which was proven in~\cite{KST}.
\begin{proposition} For any $\alpha\ge 0$
there exists some (large) constant $C=C(\alpha)$ so that the
operator $H(\tau,\sigma)$ satisfies the bounds
\begin{align}
\|H(\tau,\sigma) \|_{ L^{2,\alpha}_\rho \to L^{2,\alpha+1/2}_\rho}
&\lesssim \tau \Big(\frac{\sigma}{\tau}\Big)^{C} \label{hbd}\\
\Big\| \Big(\partial_{\tau} - 2\beta(\tau)\xi
\partial_{\xi}\Big) H(\tau,\sigma) \Big\|_{ L^{2,\alpha}_\rho \to L^{2,\alpha}_\rho}
&\lesssim  \Big(\frac{\sigma}{\tau}\Big)^{C}
\label{htaubd}\end{align} uniformly in $\sigma\ge\tau$.
\label{proph}\end{proposition}

\noindent This leads us to introduce the spaces $L^{\infty,N}
L^{2,\alpha}_\rho$ with norm
\[
\| f\|_{L^{\infty,N} L^{2,\alpha}_\rho} := \sup_{\tau\ge1} \tau^N
\|f(\tau)\|_{L^{2,\alpha}_\rho}
\]
Then the above proposition immediately allows us to draw the
following conclusions:
\begin{cor} Given $\alpha\ge0$, let $N$ be large enough. Then
\[
\| H b\|_{L^{\infty,N-2} L^{2,\alpha+1/2}_\rho}  + \Big\|
\Big(\partial_{\tau} - 2\beta(\tau)\xi
\partial_{\xi}\Big)H
b\Big\|_{L^{\infty,N-1} L^{2,\alpha}_\rho} \le C_0\,\frac{1}N
\|b\|_{L^{\infty,N} L^{2,\alpha}_\rho}
\]
with a constant $C_0$ that depends on $\alpha$ but does not depend
on $N$.
\end{cor}

\noindent The small factor $N^{-1}$ is crucial here for our argument
to work. For $H_0$ we have a stronger straightforward counterpart of
the above result:

\begin{lemma}
The operator $H_0$ satisfies the bounds
\[
\| H_0 b_d\|_{L^{\infty,N} } + \|
\partial_{\tau}H_0
b_d\|_{L^{\infty,N} } \le C_N
\|b_d\|_{L^{\infty,N}}
\]
\end{lemma}
We note that while the constant on the right cannot be small, we no
longer lose powers of $\tau$ compared to the bounds for $H$. Hence
for fixed $N$ we can choose $\tau_0$ depending on $N$ so that we
gain the smallness in the range $\tau > \tau_0$. Combining the last
two results we obtain

\begin{proposition} Given $\alpha\ge0$, let $N$ be large enough.
 Then there exists $\tau_0$ depending on $N$ so that
for $\tau > \tau_0$ we have
\[
\| \calH b\|_{L^{\infty,N-2} L^{2,\alpha+1/2}_\rho}  + \Big\|
\Big(\partial_{\tau} - 2\beta(\tau)\xi
\partial_{\xi}\Big)\calH
b\Big\|_{L^{\infty,N-1} L^{2,\alpha}_\rho} \le C_0\,\frac{1}N
\|b\|_{L^{\infty,N} L^{2,\alpha}_\rho}
\]
with a constant $C_0$ that depends on $\alpha$ but does not depend
on $N$.
\label{linmap}
\end{proposition}

On the other hand, the nonlinear operator $N_{2k-1}$
from~\eqref{final} has the  mapping properties stated in
Proposition~\ref{prop:N2k1} below. We first relate the spaces
$L^{2,\alpha}_\rho$ to the Sobolev spaces in $\R^3$.

\begin{lemma}
Let $\alpha \geq 0$. Then
\[
\| x\|_{L^{2,\alpha}_\rho} \asymp \| R^{-1} \calF^{-1}
x\|_{H^{2\alpha}(\R^3)}
\]
\end{lemma}
\begin{proof}
 For integer $k$ we have
\[
\begin{split}
\|  x\|_{L^{2,k}_\rho} & \asymp \sum_{j=0}^k \| \cL^j \calF^{-1}
x\|_{L^2} \asymp  \sum_{j=0}^k \| R^{-1}  \cL^j \calF^{-1}
x\|_{L^2(\R^3)} \\ & = \sum_{j=0}^k \| (R^{-1}  \cL R)^j R^{-1}
\calF^{-1} x\|_{L^2(\R^3)}
\end{split}
\]
But
\[
R^{-1}  \cL R = -\partial_R^2 -\frac{2}{R} \partial_R -5W^4(R)
\]
where the first two terms  can be recognized as the radial part of
the three-dimensional Laplacian. Hence we get
\[
\|  x\|_{L^{2,k}_\rho}  \asymp
 \sum_{j=0}^k \| (-\Delta - 5W^4(R))^j  R^{-1} \calF^{-1}
x\|_{L^2(\R^3)}
\]
Since $W$ is bounded together with all its derivatives, the
conclusion of the lemma follows for integer $\alpha$.

For noninteger $\alpha$ we use interpolation. First we consider the
map
\[
x \mapsto R^{-1} \calF^{-1}x
\]
and obtain the bound
\[
\| R^{-1} \calF^{-1} x\|_{H^{2\alpha}(\R^3)} \lesssim \|
x\|_{L^{2,\alpha}_\rho}
\]
To obtain the reverse bound  we use the map
\[
u \mapsto \calF R S(u)
\]
where $S(u)$ stands for the spherical average of a function $u$ in
$\R^3$.
\end{proof}

\begin{proposition}\label{prop:N2k1} Assume that $N$ is large enough and $\frac18 \leq
  \alpha < \frac\nu4$. Then the map
\[
x \to \lambda^{-2} \calF R (N_{2k-1} (R^{-1}
\calF^{-1} x))
\]
is locally Lipschitz from $ L^{\infty,N-2} L^{2,\alpha+1/2}_\rho$ to
$L^{\infty,N} L^{2,\alpha}_\rho$.\label{propn}
\end{proposition}

\begin{proof}
Using the lemma, it remains to prove that the map
\[
v \mapsto \lambda^{-2} N_{2k-1} (v)
\]
is locally Lipschitz from $ L^{\infty,N-2} H^{2\alpha+1}$ to
$L^{\infty,N} H^{2\alpha}$. We recall that
\[
N_{2k-1}(v) =5(u_{2k-1}^4-u_0^4) \,v +
  10 u_{2k-1}^3 \,v^2 + 10 u_{2k-1}^2 \,v^3 + 5u_{2k-1} \,v^4 + v^5,
\]
see the comment following \eqref{eq:epsPDE}.  The time decay is
trivially obtained for all but the first term, for which we need an
additional step, where we pull out a factor of $u_{2k-1} - u_0$.
Using the regularity of $u_{2k-1}$ given in Theorem~\ref{thm:sec2}
we obtain
\[ u_{2k-1}- u_0 \in \frac{\lambda^{\frac12}}{(t\lambda)^2}
\IS^2(R,\cQ_{k-1}),\quad \lambda^{-2} (u_0^4 - u_{2k-1}^4) \in
\tau^{-2} S^2(R^{-2},\cQ_{k-1})
\]
This indicates that two units of decay in $\tau$ are gained. On the
level of Sobolev spaces we argue as follows: since we are working
with inhomogeneous Sobolev spaces, we can localize the above
estimate to unit cubes, as the $\ell^2$ summability for
$N_{2k-1}(v)$ with respect to unit cubes is inherited from any of
the $v$ factors. But in any unit cube $Q$ the coefficients
$u_{2k-1}$ have at most $(1-a)^{\frac{\nu+1}2-}$ singularities
(where $a=\frac{r}{t}\sim \frac{R}{\tau}$) so that we can bound them
in Sobolev spaces
\[
\| u_{2k-1}\|_{H^{1+2\alpha}(Q)} \lesssim
\frac{\lambda^{\frac12}}{\dist(Q,0)}, \qquad \| u_{2k-1} -
u_0\|_{H^{1+2\alpha}(Q)} \lesssim \frac{\lambda^{\frac12}}{\tau^{2}}
\dist(Q,0)
\]
where we used that $\alpha<\frac{\nu}{4}$.  Then it suffices to
establish the quintilinear estimate
\[
H^{2\alpha+1} \cdot H^{2\alpha+1}\cdot H^{2\alpha+1}\cdot H^{2\alpha+1}
\cdot H^{2\alpha+1} \subset H^{2\alpha}
\]
which in three space dimensions holds for $\alpha \geq \frac18$ (a
standard application of the fractional Leibnitz rule and Sobolev
imbedding, see \cite{taylor}, page~105).
\end{proof}

\section{Conclusion}

We now prove Theorem~\ref{Main}. We first construct a blow-up
solution inside the cone as follows. We begin with the approximate
solution $u_{2k-1}$ and the error $e_{2k-1}$ given by
Theorem~\ref{thm:sec2} inside the cone. We extend them outside the
cone to functions having the same size and regularity, supported in
$r < 2t$. Then the relation
  \[
  e_{2k-1} = (-\partial_t^2 + \partial_r^2 +\frac{2}r \partial_r) u_{2k-1}
  + u_{2k-1}^5
  \]
is valid only inside the cone.

Using Propositions~\ref{linmap}, \ref{propn} we iteratively find a
solution
\[
X \in L^{\infty,N-2} L^{2,\alpha+1/2}_\rho
\]
 for the equation \eqref{final1} for $t \leq t_0$ sufficiently small
and
\[
\frac18 \leq \alpha < \frac{\nu}4
\]
Then we set
\[
v = R^{-1} \calF^{-1} X \in  L^{\infty,N-2} H^{2\alpha+1}
\]
and
\[
u = u_{2k-1} + v
\]
Given the derivation of \eqref{final1}, the function $u$ solves
\[
(-\partial_t^2 + \partial_r^2 +\frac{2}r \partial_r) u
  + u^5 = (-\partial_t^2 + \partial_r^2 +\frac{2}r \partial_r) u_{2k-1}
  + u_{2k-1}^5 - e_{2k-1}
\]
which implies that the function $u$ solves the nonlinear wave
equation
\[
(-\partial_t^2 + \partial_r^2 +\frac{2}r \partial_r) u
  + u^5 = 0
\]
inside the cone.

The second part of the argument is to extend the above solution $u$
to the exterior of the cone $K=\{0<t<t_0,\; 0\le|x|\le t\}$ so that
the blow-up occurs only at the tip of the cone. For this we first
observe that the above function $u$ is close to $u_0$ inside the
cone and close to $0$ outside, namely

\[
\lim_{t \to 0} \int_{K_{t}} |\nabla(u(t) - u_0(t))|^2 + | u(t) -
u_{0}(t)|^6 \,dx = 0
\]
and
\[
\lim_{t \to 0} \int_{K_{t}^c} |\nabla u(t)|^2 + |u(t)|^6 \,dx  = 0
\]
Hence given $\delta > 0$ we can choose $t_0$ so that the two
quantities above are less than~$\delta^6$.

We let $w$ be the solution to the equation
\[
(-\partial_t^2 + \partial_r^2 +\frac{2}r \partial_r) w + w^5 = 0
\]
with initial data
\[
w(t_0) = u(t_0), \qquad w_t(t_0) = u_t(t_0)
\]
Due to the finite speed of propagation we conclude that
$w = u$ inside the cone. To conclude the proof of the theorem we
will show that $w$ cannot blow up outside the cone before time $0$.
For this it suffices to prove that the energy outside the cone stays
small,
\begin{equation}
\int_{K_{t}^c} |\nabla w(t)|^2 + |w(t)|^6 \, dx  \lesssim \delta
\label{outenergy}\end{equation} see~\cite{SS}, \cite{Sogge}.

This is proved using energy conservation. The energy
$\calE(w(t))$ is conserved in time. At time $t_0$ we have
\[
\calE(w(t_0)) = \calE(u(t_0)) =  \calE(W) + O(\delta)
\]
hence at time $t$ we must have a similar relation.
But the energy inside the cone is already close to this,
so we obtain
\[
\left | \int_{K_t^c} \frac12(u_t^2+|\nabla u|)^2 - \frac{|u|^6}6
\,dx \right| \lesssim \delta
\]
On the other hand, we have the Sobolev inequality
\[
\int_{K_t^c} |u|^6 dx \lesssim \left(\int_{K_t^c}|\nabla u| ^2 \,dx
\right)^3
\]
with a universal constant (independent of~$t$).  Combining the two
inequalities above we see that for each $t$ there are two
possibilities. Either we have
\[
 \int_{K_t^c} \frac12(u_t^2+|\nabla u|)^2 \, dx\lesssim \delta
\]
or
\[
 \int_{K_t^c} \frac12(u_t^2+|\nabla u|)^2 \, dx \gtrsim 1
\]
The first alternative holds at $t = t_0$. Then a continuity argument
shows that it must hold at all $t$, since the above integral is a
continuous function of $t$ for as long as it stays small.

\end{document}